\documentclass[runningheads]{llncs}
\usepackage{graphicx}
\usepackage{latexsym}
\usepackage{amssymb}
\usepackage{amsmath}
\usepackage{mathrsfs}
\usepackage{enumerate}
\usepackage{booktabs}
\usepackage{bm}
\usepackage[all]{xy}
\usepackage{color}

\DeclareFontEncoding{LS1}{}{}
\DeclareFontSubstitution{LS1}{stix2}{m}{n}

\newcommand{\blackdiamond}{{\text{\usefont{LS1}{stix2bb}{m}{it}\symbol{"E0}}}}

\newcommand{\form}{\text{Form}}
\newcommand{\fin}{\text{fin}}

\newcommand{\F}{\mathfrak{F}}
\newcommand{\M}{\mathfrak{M}}

\newcommand{\vsim}{\mid\!\sim}

\begin{document}
\title{Fundamental Propositional Logic with Strict Implication}
\titlerunning{Fundamental Propositional Logic with Strict Implication}
\author{Zhicheng Chen}
\authorrunning{Z. Chen}
\institute{Institute of Philosophy, Chinese Academy of Sciences, Beijing, China\\School of Humanities, University of Chinese Academy of Sciences, Beijing, China}
\maketitle     

\begin{abstract}

Fundamental logic was introduced by Wesley Holliday \cite{holliday_fundamental_2023} to unify intuitionistic logic and quantum logic from a proof-theoretic perspective, capturing the logic determined solely by the introduction and elimination rules of connectives $\neg$, $\wedge$, $\vee$. 
This paper incorporates strict implication---standard in intuitionistic logic and a significant candidate for quantum logic---into the framework of fundamental propositional logic. We demonstrate that, unlike the original language, the presence of strict implication causes the semantic consequence relations over pseudo-reflexive pseudo-symmetric frames and reflexive pseudo-symmetric frames to diverge. Consequently, we provide separate axiomatizations for these two logics in the language $\{\perp, \wedge, \vee, \rightarrow\}$. Soundness and completeness theorems are established for both systems.


\end{abstract}

\section{Introduction}

The history of non-classical logic has witnessed the independent emergence of two prominent strands: intuitionistic logic (\textbf{IL}) and quantum logic (\textbf{QL}). 
On one side, \textbf{IL} crystallized from the constructivist critique of classical mathematics initiated by Brouwer and formalized by Heyting \cite{Brouwer1907-ODGDW,heyting_formalen_1930}, rejecting the Law of Excluded Middle and Double Negation Elimination in favor of a semantics grounded in verification and constructive proof \cite{dummett2000elements}. On the other, the tradition of \textbf{QL} emerged from Birkhoff and von Neumann's analysis of the empirical necessities of quantum mechanics \cite{birkhoff_logic_1936}. 
For decades, these systems were viewed as disparate, if not antagonistic. 
The intersection of these logics---the ``common ground''---remained under-explored until recent theoretical developments identified certain substrate logics capable of unifying them.

A significant step toward bridging this divide was the development of a unified relational semantics. As observed by Dalla Chiara and Giuntini \cite[pp.~139--140]{gabbay_quantum_2002}, it is possible to construct a general Kripkean framework where valuations are restricted to satisfy a closure property—specifically, \begin{align}
\Box_{\F} \blackdiamond_{\F} V(p) \subseteq V(p) \textit{, for each proposition letter } p \tag{*} \label{*}
\end{align}
where $\blackdiamond_{\F} X=\{w\in \F\mid \textit{there is a predecessor of } w \textit{ in }X\}$ and $\Box_{\F} X=\{w\in \F \mid $ every successor of $w$ is in $X\}$, for each  each set $X$ of worlds in the given frame $\F$. Let us call the Kripke models satisfying \eqref{*} as ``$\Box\blackdiamond$-models". 
Within this semantic architecture, the distinction between these major non-classical traditions reduces to simple frame conditions: requiring the accessibility relation to be reflexive and transitive yields the Kripke semantics for \textbf{IL} \cite{Kripke1963-SemanticalIntuitionistic}, whereas requiring reflexivity and symmetry yields the semantics for orthologic (\textbf{OL}), the minimal quantum logic \cite{goldblatt_semantic_1974}. Based on this observation, several weak logics have been identified that are strictly weaker than both \textbf{IL} and \textbf{OL} \cite{chen2025unified,zhong_general_2021,zhong_propositional_2025}.

More recently, this unification has been approached from a distinct, proof-theoretic perspective by Wesley Holliday. In his introduction of ``fundamental logic" \cite{holliday_fundamental_2023}, Holliday aims to capture the intrinsic logic of the connectives $\wedge$, $\vee$, and $\neg$ as determined solely by their introduction and elimination rules in Fitch-style natural deduction. By focusing on these inferential rules, fundamental logic identifies an elegant base system that unifies \textbf{OL} and the $\{\wedge, \vee, \neg\}$-fragment of \textbf{IL}. 

Algebraically, fundamental propositional logic (\textbf{FPL}) is characterized by bounded lattices equipped with weak pseudocomplementations---a structure previously studied in different contexts. The term itself is adopted by Holliday from the work of Dzik et al. \cite{dzik2006relational,dzik2006relational:survey}. Interestingly, lattices of this type appear even earlier in the literature: they were investigated by J\"{u}rgen Schulte M\"{o}nting \cite{schulte_monting_cut_1981} under the name ``almost orthocomplemented lattices''. M\"{o}nting also provided a cut-eliminable sequent calculus for the corresponding logic (i.e., \textbf{FPL}), which can be viewed as an intersection of the sequent calculi for \textbf{OL} and \textbf{IL}. Unaware of M\"{o}nting's earlier work, Juan Aguilera et al. independently rediscovered this sequent system and used it to establish the decidability of fundamental logic \cite{aguilera2025fundamental}.

In terms of relational semantics, fundamental logic also converges with the insights of Dalla Chiara and Giuntini, as it relies on the very same class of structures---the $\Box\blackdiamond$-models.
Specifically, based on his study of representations of lattices with negations in \cite{Holliday2022-Compatibility-and-Accessibility}, Holliday shows that \textbf{FPL} is sound with respect to the class of pseudo-reflexive and pseudo-symmetric $\Box\blackdiamond$-models, and complete with respect to the subclass of reflexive and pseudo-symmetric $\Box\blackdiamond$-models \cite{holliday_fundamental_2023}. 


However, the landscape changes when we consider implication. 
In \cite{holliday_fundamental_2023}, Holliday  examines various ways to extend the relational semantics of \textbf{FPL} to treat conditionals. 
In particular, he argues that the following kind, which he calls ``preconditional'', is better than other options from a natural language point of view.
$$\lVert\varphi\hookrightarrow\psi\rVert_{\M}=\Box_{\F} (-\lVert\varphi\rVert_{\M}\cup\blackdiamond_{\F} (\lVert\varphi\rVert_{\M}\cap\lVert\psi\rVert_{\M}))$$
where $\lVert\varphi\rVert_{\M}$ denotes the truth set of $\varphi$ in the model $\M$.
Moreover, from a technical perspective, this kind of conditional encompasses both the strict implication in \textbf{IL} and the ``Sasaki hook'' over \textbf{OL} (see, e.g., \cite[pp.~147--149]{gabbay_quantum_2002}). For example, one can check that $\lVert\varphi\hookrightarrow\psi\rVert_{\M}=\Box_{\F} (-\lVert\varphi\rVert_{\M}\cup\lVert\psi\rVert_{\M})$ always holds in reflexive transitive $\Box\blackdiamond$-models. 
The axioms of the minimal ``preconditional logic'' are given in \cite{holliday_fundamental_2023} and a more detailed axiomatic study of preconditionals can be found in \cite{holliday2024preconditionals}.

Nevertheless, setting aside natural language considerations, it is technically very natural to consider extending \textbf{FPL} with strict implication, as it is semantically transparent and is the standard implication of \textbf{IL}. Furthermore, strict implication has been advocated as a natural candidate implication for quantum logic (\cite[pp.~150--151]{gabbay_quantum_2002}), and recent studies suggest it holds distinct advantages over other candidate implications in the context of \textbf{OL} \cite{kawano_sequent_2021}.

In this paper, we incorporate strict implication into \textbf{FPL}. We demonstrate that, unlike in the original language, the semantic consequence relations over the class of pseudo-reflexive pseudo-symmetric $\Box\blackdiamond$-models (denoted as $\mathcal{D}_1$) and the class of reflexive pseudo-symmetric $\Box\blackdiamond$-models (denoted as $\mathcal{D}_2$) diverge when strict implication is present. Consequently, we provide separate axiomatizations for the logics of $\mathcal{D}_1$ and $\mathcal{D}_2$ in the language $\{\bot, \wedge, \vee, \to\}$ (negation is defined using $\to$ and $\bot$), proving the corresponding soundness and completeness theorems.
Our results not only settle the axiomatization problem for \textbf{FPL} with strict implication but also deepen the understanding of how implication interacts with non-distributive and non-classical structures. 

The paper is structured as follows. Section 2 lays out the formal semantics, introducing the logics of $\mathcal{D}_1$ and $\mathcal{D}_2$, and the basic syntactic rules common to both systems. Section 3 presents the axiomatizations $\,\vdash_1$ and $\,\vdash_2$, extending the basic rules with conditions for pseudo-reflexivity, pseudo-symmetry, and reflexivity. A key technical device is the ``i-formula'' calculus, which---together with specific fixpoint rules---bridges the gap between the syntax and the $\Box\blackdiamond$-closure semantics. Section 4 proves soundness and completeness via canonical models built from pairs of sets. An appendix provides a completeness proof for the minimal $\Box\blackdiamond$-logic $\vdash_{\mathbf{K}}$.





\section{Preparations}
\begin{definition}[Formulas]
Let $PL=\{p_0,p_1,...\}$ be a countable set of propositional letters. 
The formulas are built from the logical constant $\bot$ and propositional letters, using connectives $\land,\lor$ and $\to$. Denote by $Form$ the set of all formulas. 
We will use Greek letters such as $\alpha,\beta,\varphi,\psi$ to denote arbitrary formulas, and capital Greek letters such as $\Gamma,\Delta,\Phi$ to denote arbitrary sets of formulas.

For convenience, we adopt the following syntactic conventions:
\begin{enumerate}
    \item[] $\neg \alpha := \alpha\to\bot$,
    \item[] $\top := \bot\to\bot$.
\end{enumerate}    
\end{definition}
And the order of precedence among connectives is $\neg>\land=\lor>\,\to$.

\subsection{Formal semantics}

\begin{definition}\label{definition of basic frame operators}
Let \( \mathfrak{F} = \langle W, R \rangle \) be a Kripke frame.  
We define two frame operators on $\mathfrak{F}$: for each \( X\subseteq W \),
 \begin{itemize}
     \item[] $\blackdiamond_{\mathfrak{F}}X=\{w\in W\mid \textit{there is an } x \in X \textit{ with }xRw\}$,
     \item[] $\Box_{\mathfrak{F}}X = \{ w \in W \mid$ for all $v$ with $wRv$, $v\in X \}$.
 \end{itemize}
\end{definition}

It is easy to show that $\blackdiamond_{\mathfrak{F}}X\subseteq Y \Leftrightarrow X\subseteq\Box_{\mathfrak{F}}Y$. Thus, $\blackdiamond_{\mathfrak{F}}$ and $\Box_{\mathfrak{F}}$ form a Galois connection between $\langle \wp(W),\subseteq\rangle$ and itself. It follows that $\Box_{\mathfrak{F}}\blackdiamond_{\mathfrak{F}}$ is a closure operator on $\wp(W)$, and $FP_{\mathfrak{F}}=\{\Box_{\mathfrak{F}}X \mid X\subseteq W\}$ is the set of fixpoints of $\Box_{\mathfrak{F}}\blackdiamond_{\mathfrak{F}}$ (cf. \cite{davey2002introduction} (pp.155--160)). Then, by the property of closure operators (see, e.g., \cite{davey2002introduction} (Thm. 7.3)), we have:

\begin{proposition}\label{FP complete lattice}
$FP_{\mathfrak{F}}$ ordered by $\subseteq$ forms a complete lattice with the least element $\Box_{\mathfrak{F}}\emptyset$ and the greatest element $W$, and  
        \[\bigwedge_{i\in I} X_i = \bigcap_{i\in I} X_i \text{ and } \bigvee_{i\in I} X_i = \Box_{\mathfrak{F}}\blackdiamond_{\mathfrak{F}}(\bigcup_{i\in I} X_i). \].     
\end{proposition}

\begin{definition}[Connectives as frame operators]\label{def: connectives as frame operators}
Given a frame \( \mathfrak{F}= \langle W,R\rangle\), let us define the following operators on $FP_{\mathfrak{F}}$: for any $X,Y\subseteq W$,
\begin{itemize}
    \item[] $\bot_{\mathfrak{F}} = \Box_{\mathfrak{F}}\emptyset$,
    \item[] $X \to_{\mathfrak{F}} Y = \Box_{\mathfrak{F}} (-_{\mathfrak{F}}X \cup Y)$,
    \item[] $X \lor_{\mathfrak{F}} Y = \Box_{\mathfrak{F}}\blackdiamond_{\mathfrak{F}}(X \cup Y)$,
\end{itemize}
where $-_{\mathfrak{F}}X = W\setminus X$. 
\end{definition}

\begin{definition}[Interpretation of formulas]\label{def: interpretation of formulas}
Given a Kripke model \( \mathfrak{M}\), we define the truth sets of formulas in $\mathfrak{M}$ recursively as follows:
\begin{itemize}
    \item $\| p \|_{\mathfrak{M}} = V(p)$, for any $p \in PL$;
    \item $\| \bot \|_{\mathfrak{M}}= \bot_{\mathfrak{F}}$;
    \item $\| \alpha \land \beta \|_{\mathfrak{M}} = \| \alpha \|_{\mathfrak{M}} \cap \| \beta \|_{\mathfrak{M}}$;
    \item $\|\alpha \to \beta\|_{\mathfrak{M}} = \| \alpha \|_{\mathfrak{M}}\to_{\mathfrak{F}} \| \beta \|_{\mathfrak{M}}$;
    \item $\|\alpha \lor \beta\|_{\mathfrak{M}} = \| \alpha \|_{\mathfrak{M}}\lor_{\mathfrak{F}} \| \beta \|_{\mathfrak{M}}$.
\end{itemize}

Then for any \(w\in\mathfrak{M}\) and formula \( \varphi \), let
\[
\mathfrak{M}, w \vDash \varphi \iff w \in \| \varphi \|_{\mathfrak{M}}.
\]
\end{definition}

For simplicity, and where no ambiguity arises, we will omit the subscript from the notation “$\Vert\cdot\Vert$” 
(except in the formal statements of the definitions, lemmas, and theorems). 

\begin{definition}[$\Box\blackdiamond$-model]\label{def: DG-model and interpretation of formulas}
  A Kripke model $\mathfrak{M}=\langle \mathfrak{F},V\rangle$ is called a $\Box\blackdiamond$-model if $\,V$ maps from $PL$ to $FP_{\mathfrak{F}}$. 
  Denote by $\mathcal{D}_{\textbf{$\Box\blackdiamond$}}$ the class of $\Box\blackdiamond$-models.
\end{definition}

By an easy induction, we can show that:
\begin{lemma}\label{formula.fixedpoint}
Let $\mathfrak{M}=\langle \mathfrak{F},V\rangle$ be a $\Box\blackdiamond$-model. For each formula $\varphi$, $\Vert \varphi \Vert_{\mathfrak{M}}\in FP_{\mathfrak{F}}$.
\end{lemma}

\begin{definition}\label{def: PsRe, PsSy}
Let us define the following notions. (Note that $\Diamond_{\mathfrak{F}}$ is the dual of $\Box_{\mathfrak{F}}$, i.e.,  $-_{\mathfrak{F}}\Box_{\mathfrak{F}}-_{\mathfrak{F}}$.)
    \begin{itemize}
    \item A frame $\mathfrak{F}$ is pseudo-reflexive iff for each $w$ in $\mathfrak{F}$, $w\in \Box_{\mathfrak{F}}\emptyset\cup\Diamond_{\mathfrak{F}}\Box_{\mathfrak{F}}\blackdiamond_{\mathfrak{F}}\{w\}$. 
    \item A frame $\mathfrak{F}$ is pseudo-symmetric iff for each $w$ in $\mathfrak{F}$, $w\in \Box_{\mathfrak{F}}\Diamond_{\mathfrak{F}}\Box_{\mathfrak{F}}\blackdiamond_{\mathfrak{F}}\{w\}$.
    \end{itemize}
And a model is pseudo-reflexive/pseudo-symmetric when its underlying frame is.
\end{definition}

\begin{definition}[Semantic consequence]
Let $\mathcal{D}_{1}$ be the class of  pseudo-reflexive pseudo-symmetric $\Box\blackdiamond$-models, and  $\mathcal{D}_{2}$ the class of reflexive  pseudo-symmetric $\Box\blackdiamond$-models.

For each $k\in\{1,2\}$, let $\mathbf{F}_k$ denote the logic of $\mathcal{D}_k$,  with $\,\vDash_{k}$ being the semantic consequence relation  --- for each set \( \Gamma \) of formulas and each formula \( \varphi \), 
\[
\Gamma \vDash_{k} \varphi \iff \text{for every } \mathfrak{M} \in \mathcal{D}_k \text{ and } w \in \mathfrak{M},\ (\mathfrak{M}, w \vDash \Gamma \Rightarrow \mathfrak{M}, w \vDash \varphi).\]
\end{definition}

In the original language $\{\bot,\land,\lor,\neg\}$ of \textbf{FPL}, the semantic consequence relations over $\mathcal{D}_1$ and $\mathcal{D}_2$ are identical. However, it is not the case in our language $\{\bot,\land,\lor,\to\}$, as shown in the following lemma.
\begin{lemma}
    $\,\vDash_{1}\subsetneq \,\vDash_{2}$.
\end{lemma}
\begin{proof}
    $\,\vDash_{1}\subseteq \,\vDash_{2}$ follows easily from that fact that $\mathcal{D}_{2}\subseteq \mathcal{D}_{1}$. It remains to prove $\,\vDash_{1}\neq \,\vDash_{2}$. On one hand, note that $\{\alpha,\alpha\to\beta\}\vDash_2\beta$ for each formula $\alpha$ and $\beta$. On the other hand, consider the following frame $\F=\langle \mathbb{N}, <\rangle$. It is easy to check that $\F$ is  pseudo-reflexive and pseudo-symmetric, and that for each $n\in\mathbb{N}$, $\{j\mid n\leq j\}\in FP_{\mathfrak{F}}$. Take two propositional letters $p$ and $q$. Define a valuation function $V$: $V(q)=\mathbb{N}\setminus\{0\}$ while $V(r)=\mathbb{N}$ for any other propositional letter $r$. Let $\M =\langle \F,V\rangle$, then we can show that $\M\in\mathcal{D}_1$, $\M ,0\vDash p$ , $\M ,0\vDash p \to q$, and $\M ,0\nvDash q$. So, $\{ p ,  p \to q\}\nvDash_1 q$. As a consequence, $\,\vDash_{1}\neq \,\vDash_{2}$.
\end{proof}

We end this section with a discussion of our definition of negation, namely ``$\neg\alpha:=\alpha\to\bot$''. Let us call a frame $\F$ ``successor-serial'' if, in $\F$, every successor is not a dead end. The next lemma  shows that, semantically speaking, the negation defined in this way is a ``strict negation'' on successor-serial frames. 
\begin{lemma}\label{semantics: neg = to bot}
     If a frame \( \mathfrak{F}= \langle W,R\rangle\) is successor-serial, then for each $X\subseteq W$,
    $$X \to_{\mathfrak{F}} \bot_{\mathfrak{F}} = \Box_{\mathfrak{F}} -_{\mathfrak{F}}X.$$
\end{lemma}

Furthermore, one can easily see that pseudo-symmetric frames are successor-serial. The same holds for reflexive frames.


\subsection{Basic syntactic rules}
Generally, our definition of a syntactic consequence relation proceeds in two steps. First, we define a binary relation $\,\vdash\,\subseteq Form \times Form$ via specified syntactic rules. 
We then extend it to sets of formulas in the standard way: for each $\Gamma\subseteq Form$ and $\varphi\in Form$, 
\[\Gamma\vdash^+ \alpha \iff \text{there exists finite } \Gamma_0\subseteq\Gamma \text{ such that }\bigwedge\Gamma_0\vdash\alpha\]  with the convention that $\bigwedge\emptyset=\top$. For simplicity, hereater we shall abuse notation by writing $\,\vdash$ for  $\,\vdash^+$.


To begin with, let us look at the following properties for a  relation $\, \vdash\ \subseteq Form \times Form$, which are called the ``basic rules''. (For brevity, universal quantifiers over $\Gamma $, $\alpha$, etc are left tacit in rule statements; the subscript of $\bigwedge$ and $\bigvee$ is simplied when clear from context; and $\,\vdash\alpha$ abbreviates ``$\chi\vdash\alpha$ for every formula $\chi$''.)
\begin{definition}[Basic rules]\label{basic rule}
    \begin{itemize}
    \item[] (A)\quad \(  \alpha  \vdash \alpha \)
    \item[] (Cut)\quad \(   \alpha  \vdash \beta \text{ and } \beta  \vdash\gamma\ \Rightarrow\ \alpha  \vdash \gamma \)
    
    \item[] (\(\bot\))\quad \( \bot \vdash \alpha \)

    \item[] (\(\land\)R)\quad \( \chi\vdash\alpha  \) and \( \chi\vdash\beta\   \Rightarrow\  \chi\vdash\alpha\land\beta \)
    \item[] (\(\land\)L)\quad \( \alpha \land \beta \vdash \alpha \) and \( \alpha \land \beta \vdash \beta  \)

    \item[] (\(\lor\)R)\quad \(  \alpha \vdash \alpha \lor \beta \) and \( \beta  \vdash \alpha \lor \beta \)
    \item[] (\(\lor\)L)\quad \( \alpha \vdash \chi \) and \( \beta \vdash \chi\   \Rightarrow\  \alpha\lor \beta \vdash \chi \)
    
    \item[] (DT${}_0$)\quad $\alpha \vdash \beta\  \Rightarrow \,\  \vdash\alpha\to\beta$ 
    \item[] ($\to\land$)\quad $(\alpha\to\beta)\land (\alpha\to\gamma)\,\vdash\alpha\to \beta\land\gamma$ 
    \item[] ($\to$tr)\quad $(\alpha\to\beta)\land (\beta\to\gamma)\,\vdash\alpha\to \gamma$ 

    \item[] ($\rightarrow$-$\lor$.s) \quad 
    \(\bigvee_{1\leq j\leq n}(\alpha_j \rightarrow \beta_j)  \land \bigwedge_j(\psi_j \land \beta_j \rightarrow \chi ) \vdash \bigwedge_j(\psi_j \land \alpha_j) \rightarrow \chi\) (for each $n\geq 1$)
\end{itemize}
\end{definition}

Specifically, the rules mentioning $\to$ are called ``basic $\to$-rules''.

It is readily verified that:
\begin{lemma}\label{derived rules of basic rules}
    If $\,\vdash$ satisfies basic rules, then the following holds:
    \begin{itemize}
        \item[] (Mon)\quad \( \Gamma \vdash \varphi \text{ and } \Gamma\subseteq\Delta  \Rightarrow \Delta \vdash \varphi \)
        \item[]   ($\top$)\qquad $\vdash \top$

        \item[] ($\rightarrow re$)\quad $\vdash \varphi \rightarrow \varphi$
        \item[] ($\land \rightarrow \land$)\quad $\bigwedge_{1\leq j\leq n}(\alpha_j \rightarrow \beta_j) \vdash \bigwedge_j\alpha_j  \rightarrow \bigwedge_j\beta_j$
        \item[] ($\neg$ antitone)\quad $\alpha\vdash\beta \Rightarrow \neg\beta\vdash \neg\alpha$
    \end{itemize}
\end{lemma}


\begin{definition}[$\,\vDash_{\mathfrak{M}}$]
    Given a Kripke model $\mathfrak{M}$, we define a relation $\,\vDash_{\mathfrak{M}}\subseteq \wp(Form)\times Form$: for any $\Gamma\subseteq Form$ and $\varphi\in Form$, 
    \begin{align*}
       \Gamma\vDash_{\mathfrak{M}}\varphi \iff \textit{for every }w\in\mathfrak{M},\ \mathfrak{M},w\vDash\Gamma \textit{ implies } \mathfrak{M},w\vDash\varphi\ .
    \end{align*}
\end{definition}

\begin{definition}[Correctness]\label{def: correctness of rule}
    A sequent rule ($X$) is correct at a Kripke model $\mathfrak{M}$ if and only if $\ \vDash_{\mathfrak{M}}$ satisfies ($X$).   
\end{definition}

\begin{lemma}\label{correctness part.1.1}
~
   \begin{enumerate}
       \item All the basic rules, except ($\bot$) and ($\lor$L), are correct at any Kripke model. 
       
       \item ($\bot$) and ($\lor$L) are correct at any $\Box\blackdiamond$-model.
       
   \end{enumerate}
\end{lemma}

\begin{proof}
~
    \begin{enumerate}
       \item Let us show, for example, the correctness of (DT${}_0$) and ($\rightarrow$-$\lor$.s) at an arbitrary Kripke model $\mathfrak{M}$.
       \begin{itemize}
           \item Suppose $\alpha \vDash_{\mathfrak{M}} \beta$. For each $w\in\mathfrak{M}$ and each successor $v$ of $w$, since $\alpha \vDash_{\mathfrak{M}} \beta$, we know that $\mathfrak{M},v\vDash\alpha$ implies $\mathfrak{M},v\vDash\beta$. So $\,\vDash_{\mathfrak{M}}\alpha\to\beta$.
           
           \item  Note that $\Vert \bigvee(\alpha_j \rightarrow \beta_j) \Vert = \Box_{\F}\blackdiamond_{\F} \bigcup \Box_{\F}(-\Vert\alpha_j\Vert\cup\Vert\beta_{j}\Vert) = \Box_{\F}\bigcup\blackdiamond_{\F} \Box_{\F}(-$ $\Vert\alpha_j\Vert\cup\Vert\beta_{j}\Vert)  \subseteq \Box_{\F}\bigcup(-\Vert\alpha_j\Vert\cup\Vert\beta_{j}\Vert )$. 
           So if \(\mathfrak{M},w\vDash\bigvee(\alpha_j \rightarrow \beta_j)\), then any successor $v$ satisfying $\bigwedge \alpha_j$ also satisfies some $\psi_j$. From this, the correctness of ($\rightarrow$-$\lor$.s) follows immediately.
           
       \end{itemize}
       
       \item This follows from the definition of $\bot_{\mathfrak{F}}$ and $\lor_{\mathfrak{F}}$, Proposition \ref{FP complete lattice}, and
       Lemma \ref{formula.fixedpoint}. 
   \end{enumerate}
\end{proof}

Let \(\, \vdash_{\mathbf{K}} \)  be the smallest  $\, \vdash$ that satisfies all the basic rules, and let $\,\vDash_{\mathbf{K}}$  be the local semantic consequence relation over the class of all $\Box\blackdiamond$-models. By lemma \ref{correctness part.1.1} and the completeness theorem in the appendix, we have:
\begin{theorem}\label{soundness and completeness of vdash_K}
    For every $\Gamma\subseteq Form$ and $\varphi\in Form$,
    \[\Gamma\,\vdash_{\mathbf{K}}\varphi \Leftrightarrow\Gamma\,\vDash_{\mathbf{K}}\varphi.\]
\end{theorem}

\section{Syntactic consequence}

Regarding pseudo-reflexivity and pseudo-symmetry, we introduce the following rules:
\begin{definition}[(Abs), ($\neg\neg$I)]\label{def: (Abs) and (DNI)}
\begin{enumerate}
    \item[] (Abs) \quad  \( \alpha\land \neg\alpha \vdash \bot \).  
    \item[] ($\neg\neg$I) \quad  \( \alpha \vdash \neg \neg \alpha \). 
\end{enumerate}    
\end{definition}

In \cite{holliday_fundamental_2023} (Prop. 4.14), Wesley proves that pseudo-reflexivity(/pseudo-symmetry) is, in a certain sense, the frame correspondence of (Abs)(/($\neg\neg$I)). However, the meaning of ``$\neg$'' in our paper is, in general, different from that in \cite{holliday_fundamental_2023}, though they coincide on successor-serial frames by Lemma \ref{semantics: neg = to bot}.
Therefore, the frame correspondence of (Abs)(/($\neg\neg$I)) need to be adjusted in our settings, as shown in the following: 
\begin{lemma}\label{correspondance of (Abs) and (DNI)}
    Let $\F$ be a frame. In each of the following pairs, (a) and (b) are equivalent:
\begin{enumerate}
    \item
      \begin{enumerate}[(a)]
          \item for each $\Box\blackdiamond$-model $\M$ based on $\F$, (Abs) is correct at $\M$.
          \item $\F$ is strongly pseudo-reflexive: for each $w$ in $\mathfrak{F}$, $w\in \Box_{\mathfrak{F}}\emptyset \cup \Diamond_{\mathfrak{F}}(-_{\mathfrak{F}}\Box_{\mathfrak{F}}\emptyset \cap \Box_{\mathfrak{F}}\blackdiamond_{\mathfrak{F}}\{w\})$. 
      \end{enumerate} 

    \item
      \begin{enumerate}[(a)]
          \item for each $\Box\blackdiamond$-model $\M$ based on $\F$, ($\neg\neg$I) is correct at $\M$.
          \item $\F$ is weakly pseudo-symmetric: for each $w$ in $\mathfrak{F}$, $w\in \Box_{\mathfrak{F}}(\Box_{\mathfrak{F}}\emptyset \cup \Diamond_{\mathfrak{F}}(-_{\mathfrak{F}}\Box_{\mathfrak{F}}\emptyset \cap \Box_{\mathfrak{F}}\blackdiamond_{\mathfrak{F}}\{w\}))$.
      \end{enumerate} 
\end{enumerate}
\end{lemma}
\begin{proof}
~
    \begin{enumerate}
        \item 
           Assume that (b) holds. Let $\M$ be a $\Box\blackdiamond$-model based on $\F$, $w\in\M$, $\M,w\vDash\alpha$, and $\M,w\vDash\neg\alpha$. We want to show $\M,w\vDash\bot$, i.e., $w\in\Vert\bot\Vert_{\M}=\Box_{\mathfrak{F}}\emptyset$. Suppose not, then by (b) $w\in\Diamond_{\mathfrak{F}}(-_{\mathfrak{F}}\Box_{\mathfrak{F}}\emptyset \cap \Box_{\mathfrak{F}}\blackdiamond_{\mathfrak{F}}\{w\})$, i.e., there exist a successor $v$ of $w$ with $v\in-_{\mathfrak{F}}\Box_{\mathfrak{F}}\emptyset \cap \Box_{\mathfrak{F}}\blackdiamond_{\mathfrak{F}}\{w\}$. Since $\Vert\alpha\Vert_{\M}$ is a $\Box_{\mathfrak{F}}\blackdiamond_{\mathfrak{F}}$-fixpoint (by Lemma \ref{formula.fixedpoint}), from $v\in \Box_{\mathfrak{F}}\blackdiamond_{\mathfrak{F}}\{w\}$ and $\M,w\vDash\alpha$ we obtain $\M,v\vDash\alpha$, which, together with $v\in-_{\mathfrak{F}}\Box_{\mathfrak{F}}\emptyset$, contradicts  $\M,w\vDash \alpha\to\bot$.
            
           Conversely, assume (b) does not hold. Then, we can find $w_0$ in $\mathfrak{F}$, such that $w_0\notin \Box_{\mathfrak{F}}\emptyset \cup \Diamond_{\mathfrak{F}}(-_{\mathfrak{F}}\Box_{\mathfrak{F}}\emptyset \cap \Box_{\mathfrak{F}}\blackdiamond_{\mathfrak{F}}\{w_0\})$. Define a valuation function $V_0$: $V_0(p_0)=\Box_{\mathfrak{F}}\blackdiamond_{\mathfrak{F}}\{w_0\}$, and $V_0(p)=\Box_{\mathfrak{F}}\emptyset$ for each $p\in PL\setminus\{p_0\}$. Then $\M_0=\langle\F,V_0\rangle$ is a $\Box\blackdiamond$-model, and $\M_0 ,w_0 \vDash p_0$. Moreover, since $w_0\notin \Box_{\mathfrak{F}}\emptyset \cup \Diamond_{\mathfrak{F}}(-_{\mathfrak{F}}\Box_{\mathfrak{F}}\emptyset \cap \Box_{\mathfrak{F}}\blackdiamond_{\mathfrak{F}}\{w_0\})$, we have $\M_0 ,w_0 \nvDash \bot$ and $\M_0 ,w_0 \vDash p_0\to\bot$. So, (Abs) is not correct at $\M_0$, as $w_0$ witnesses it.

          \item Assume that (b) holds. Let $\M$ be a $\Box\blackdiamond$-model based on $\F$,  $w\in\M$, and $\M,w\vDash\alpha$. We want to show  $\M,w\vDash\neg\neg\alpha$, i.e., $\M,w\vDash (\alpha\to\bot)\to\bot$. Suppose not, then we can find a successor $v$ of $w$ with $v\in-_{\mathfrak{F}}\Box_{\mathfrak{F}}\emptyset \cap\Vert\alpha\to\bot\Vert_{\M}$. By (b) and $v\in-_{\mathfrak{F}}\Box_{\mathfrak{F}}\emptyset$, we obtain  $v\in\Diamond_{\mathfrak{F}}(-_{\mathfrak{F}}\Box_{\mathfrak{F}}\emptyset \cap \Box_{\mathfrak{F}}\blackdiamond_{\mathfrak{F}}\{w\})$. 
            Then, since $\M,w\vDash\alpha$ and  $\Vert\alpha\Vert_{\M}$ is a $\Box_{\mathfrak{F}}\blackdiamond_{\mathfrak{F}}$-fixpoint, we have $v\in\Diamond_{\mathfrak{F}}(-_{\mathfrak{F}}\Box_{\mathfrak{F}}\emptyset \cap \Vert\alpha\Vert_{\M})$, which contradicts  $v\in\Vert\alpha\to\bot\Vert_{\M}$.
            
            Conversely, assume that (b) does not hold. Then, we can find $w_0$ in $\mathfrak{F}$ and a successor $v_0$ of $w_0$, such that $v_0\notin \Box_{\mathfrak{F}}\emptyset \cup \Diamond_{\mathfrak{F}}(-_{\mathfrak{F}}\Box_{\mathfrak{F}}\emptyset \cap \Box_{\mathfrak{F}}\blackdiamond_{\mathfrak{F}}\{w_0\})$. Define a valuation function $V_0$: $V_0(p_0)=\Box_{\mathfrak{F}}\blackdiamond_{\mathfrak{F}}\{w_0\}$, and $V_0(p)=\Box_{\mathfrak{F}}\emptyset$ for each $p\in PL\setminus\{p_0\}$. Then $\M_0=\langle\F,V_0\rangle$ is a $\Box\blackdiamond$-model, and $\M_0 ,w_0 \vDash p_0$. Moreover, since $v_0\notin \Box_{\mathfrak{F}}\emptyset \cup \Diamond_{\mathfrak{F}}(-_{\mathfrak{F}}\Box_{\mathfrak{F}}\emptyset \cap \Box_{\mathfrak{F}}\blackdiamond_{\mathfrak{F}}\{w_0\})$, we have $\M_0 ,v_0 \nvDash \bot$ and $\M_0 ,v_0 \vDash p_0\to\bot$. So, $\M_0 ,w_0 \nvDash (p_0\to\bot)\to\bot$. Thus, ($\neg\neg$I) is not correct at $\M_0$.
    \end{enumerate}
\end{proof}

Observe that:
(1) strong pseudo-reflexivity implies pseudo-reflexivity; (2) pseudo-symmetry implies weak pseudo-symmetry; (3) assuming pseudo-symmetry, pseudo-reflexivity is equivalent to strong pseudo-reflexivity. So, together with Lemma \ref{correspondance of (Abs) and (DNI)} we have:
\begin{corollary}\label{correctness part.1.2}
    ($\neg\neg$I) is correct at every pseudo-symmetric $\Box\blackdiamond$-model. And (Abs) is correct at every model in $\mathcal{D}_1$.
\end{corollary}

With regard to reflexivity, note that in reflexive models, we have $\|(\alpha_1\to\beta_1)\lor\dots\lor(\alpha_n\to\beta_n)\|$ $\subseteq$ $\blackdiamond_{\F}\|\alpha_1\to\beta_1\| \cup\dots\cup \blackdiamond_{\F}\|\alpha_n\to\beta_n\| $ $\subseteq$ $-\|\alpha_1\| \cup \|\beta_1 \| \cup\dots\cup -\|\alpha_n\| \cup\|\varphi_n \| $. So, we can introduce the following rules. (We omit the initial universal quantification ``for each $ n \in \mathbb{N}^*$'', etc.)

\begin{definition}[(Refl${}_1$), (Refl${}_2$)]\label{def: (Refl1) and (Refl2)}
    \begin{enumerate}
        \item[] (Refl${}_1$) \quad   $\psi_j \land \beta_j \vdash \chi$ for every $1 \leq j \leq n $ $\;\Rightarrow\; \bigwedge_j(\psi_j \land \alpha_j) \land \bigvee_j(\alpha_j \rightarrow \beta_j)   \vdash \chi$
        \item[] $(\text{Refl}_2)$ \quad  $\bigwedge_{1 \leq j \leq n}(\psi_j \land \beta_j\to\chi) \vdash \bigwedge_j(\psi_j \land \alpha_j) \land \bigvee_j(\alpha_j \rightarrow \beta_j)\to\chi$
    \end{enumerate}
\end{definition}

It is not hard to show that:
\begin{lemma}\label{correctness part.1.3}
  $(\text{Refl}_1)$ and $(\text{Refl}_2)$ are correct at every reflexive Kripke model.
\end{lemma}

Also notice that the rule of Modus ponens---$\alpha\land(\alpha\to\beta) \vdash \beta$ for each $\alpha$ and $\beta$---can be derived from (Refl${}_1$) by taking the special case of $n=1$. 

Now that we have introduced the basic rules as well as various rules associated with each frame condition, at first glance, it might seem we have got all we need to define $\,\vdash_k$. However, it turns out that our work on axiomatization is still missing the final piece of the ``puzzle'' --- we will introduce a relatively complex syntactic rule for each $\,\vdash_k$. This rule is related to the way in which the fact that ``the truth set of every formula in a $\Box\blackdiamond$-model is a $\Box\blackdiamond$-fixpoint'' behaves in classes of $\Box\blackdiamond$-models
with sufficiently strong properties (at least being successor-serial). To facilitate understanding, let us first look at two examples.




\begin{example}\label{example 1}
Consider $\Gamma_1=\{(\top\rightarrow\varphi_1 ) \lor (\top\rightarrow\varphi_2 ) \lor \gamma, \ (\varphi_1 \rightarrow \bot) \lor \gamma, \  (\varphi_2 \rightarrow \bot) \lor \gamma\}$. For any successor-serial $\Box\blackdiamond$-model $\mathfrak{M}$, 

 \noindent \quad $\|\bigwedge\Gamma_1\|$\\
=\ $\Box_{\F}\blackdiamond_{\F}(\|\top\rightarrow\varphi_1\| \cup \|\top\rightarrow\varphi_2 \| \cup \|\gamma\|)\ \cap\ \Box_{\F}\blackdiamond_{\F}(\|\varphi_1 \rightarrow \bot\| \cup \|\gamma\|)\ \cap\ \Box_{\F}\blackdiamond_{\F}(\|\varphi_2 \rightarrow \bot\| \cup \|\gamma\|)$\\
=\ $\Box_{\F}(\ ((\blackdiamond_{\F}\|\top\rightarrow\varphi_1 \| \cup \blackdiamond_{\F}\|\top\rightarrow\varphi_2 \|) \ \cap\ \blackdiamond_{\F}\|\varphi_1 \rightarrow \bot\| \ \cap\ \blackdiamond_{\F}\|\varphi_2 \rightarrow \bot\|)\ \cup\ \blackdiamond_{\F}\|\gamma\|\ )$ \quad   \hfill (since $\blackdiamond_{\F}$ is distributive over $\cup$ and $\Box_{\F}$ distributive over $\cap$)\\
$\subseteq$\ $\Box_{\F}(\ ((\|\varphi_1 \| \cup \|\varphi_2 \|) \ \cap\ -\|\varphi_1 \| \ \cap\ -\|\varphi_2 \|)\ \cup\ \blackdiamond_{\F}\|\gamma\|\ )$ \quad   \hfill (by Lemma \ref{semantics: neg = to bot} and the fact that $\blackdiamond_{\F}\Box_{\F} X\subseteq X$ for every $X\subseteq W$)\\
=\ $\Box_{\F}\blackdiamond_{\F}\|\gamma\|$\\
=\ $\|\gamma\|$ \qquad  (by Lemma \ref{formula.fixedpoint})\\
As a corollary,  we have $\Gamma_1 \vDash_{1} \gamma$.    
\end{example}

\begin{example}\label{example 2}
Note that in reflexive models  $\|\bigvee(\alpha_j\to\beta_j)\lor\gamma\|$ $\subseteq$ $\bigcup\blackdiamond_{\F}\|\alpha_j\to\beta_j\| \cup \blackdiamond_{\F}\|\gamma\|$ $\subseteq$ $\bigcup(-\|\alpha_j\| \cup \|\beta_j \|) \cup \blackdiamond_{\F}\|\gamma\|$. Also note that reflexive models are successor-serial. 
Consider $\Gamma_2=\{(\top\rightarrow\alpha_1) \lor \gamma, (\top\rightarrow\alpha_2) \lor \gamma, \ (\top\rightarrow((\alpha_1\to\beta_1)\lor(\alpha_2\to\beta_2)\lor\gamma)) \lor \gamma, \ (\beta_1 \rightarrow \bot) \lor \gamma, \  (\beta_2 \rightarrow \bot) \lor \gamma\}$.  Analogous to the proof of Example \ref{example 1}, one can show that for each reflexive $\Box\blackdiamond$-model, $\Vert \bigwedge \Gamma_2 \Vert \subseteq \Box_{\F}((\|\alpha_1 \| \cap \|\alpha_2 \| \cap  (-\|\alpha_1\| \cup -\|\alpha_2\| \cup \|\beta_1 \| \cup \|\beta_2 \|)  \cap -\|\beta_1 \| \cap -\|\beta_2 \|) \cup \blackdiamond_{\F}\|\gamma\|) = \Box_{\F}\blackdiamond_{\F}\|\gamma\|\subseteq\|\gamma\|$. As a corollary,  $\Gamma_2 \vDash_{2} \gamma$. 
\end{example}

In general, Example \ref{example 1} provides the following inspiration. Over  successor-serial $\Box\blackdiamond$-models, 
if we fix a formula $\gamma$, then, for every formula of the form $(\varphi_1\to\psi_1)\lor...\lor(\varphi_n\to\psi_n)\lor\gamma$, we can write a classical implication ``$\{\varphi_1, ..., \varphi_n\}\sqsupset\{\psi_1, ..., \psi_n\}$'', namely, the conjunction of $\varphi_i $ classically implying the classical disjunction of $\psi_i$. If, from a certain set of these classical implications, we can deduce falsity, then $\gamma$ should be the consequence of the set of the corresponding original formulas.
Furthermore, Example \ref{example 2} reveals that, over reflexive $\Box\blackdiamond$-models, an adequate calculus for these classical implications---relative to the given formula $\gamma$---should contain the following axiom:
$$ \{\varphi_1, \dots, \varphi_n, \bigvee\nolimits_j(\varphi_j \rightarrow \psi_j) \lor \gamma\} \sqsupset \{\psi_1, \dots, \psi_n\}$$
which also evidences that, under this circumstance,  it is not just the classical propositional calculus that we need here on these classical implications.


\begin{definition}
    Let $\form^i = \{\Delta \sqsupset \Theta \mid \Delta, \Theta \in \mathscr{P}_{\fin}^*(\form)\}$, where $\mathscr{P}_{\fin}^*(\form)$ denotes the family of non-empty  finite subsets of $Form$. The elements of $\form^i$ are called ``i-formulas''. From now on, we will use symbols with superscript ``i'', such as $\alpha^{i}$ and $\Gamma^{i}$, for i-formulas and sets of i-formulas.
\end{definition}

To define $\vdash_1$, it turns out that the following deductive system on i-formulas, which is essentially a subsystem of the classical sequent calculus (viewing each i-formula as a sequent), is sufficient.
\begin{definition}[$\,\Vdash_1$]\label{definition of Vdash_1}
     Let $\,\Vdash_1\, \subseteq \wp(\form^i) \times \form^i$ be the least $\,\Vdash$ that satisfies the following rules, which are called the ``basic i-rules'':
    \begin{enumerate}
        \item[] (A) \quad $\Gamma^i \cup \{\alpha^i\} \Vdash \alpha^i$ 
        \item[] (Cut)\quad $\Gamma^i \Vdash \alpha^i$ and $\Phi^i \cup \{a^i\} \Vdash \beta^i \,\Rightarrow\, \Gamma^i \cup \Phi^i \Vdash \beta^i$ 
        \item[] (i-A) \quad $\Delta \cap \Theta \neq \emptyset \ \Rightarrow\ \ \Vdash \Delta \sqsupset \Theta$
        \item[] (i-Cut) \quad $\{\Delta_1 \sqsupset \Theta_1 \cup \{\varphi\},\ \Delta_2 \cup \{\varphi\} \sqsupset \Theta_2\} \Vdash \Delta_1 \cup \Delta_2 \sqsupset \Theta_1 \cup \Theta_2$
        \item[] (i-$\land$L) \quad $\Vdash \{\varphi\land \psi\} \sqsupset \{\varphi \land \psi\}$, and $\Vdash \{\varphi\} \sqsupset \{ \psi\}$
        \item[] (i-$\land$R) \quad $\Vdash \{\varphi, \psi\} \sqsupset \{\varphi \land \psi\}$
    \end{enumerate}
    
\end{definition}

It is readily verified that:
\begin{lemma}\label{properties of vdash^i}
If $\,\Vdash$ satisfies the basic i-rules, then it enjoys the following properties:
\begin{enumerate}
    \item[] $(Mon)$\quad \( \Gamma^{i} \Vdash \alpha^{i} \text{ and } \Gamma^{i}\subseteq\Phi^{i}  \Rightarrow \Phi^{i} \Vdash \alpha^{i} \);
    \item[] $(i$-$Mon)$\quad $\Delta_1 \subseteq \Delta_2$ and $\Theta_1 \subseteq \Theta_2 \Rightarrow \Delta_1 \sqsupset \Theta_1 \Vdash \Delta_2 \sqsupset \Theta_2$;
    \item[] $(i$-$Cut.1)$\quad $\{\Delta_1 \sqsupset \Theta_1 \cup \Theta\} \cup \{\Delta_2 \cup \{p\} \sqsupset \Theta_2 \mid p \in \Theta\} \Vdash \Delta_1 \cup \Delta_2 \sqsupset \Theta_1 \cup \Theta_2$;
    \item[] $(i$-$Cut.2)$\quad $\{\Delta_1 \sqsupset \Theta_1 \cup \{\varphi\} \mid \varphi \in \Delta\} \cup \{\Delta_2 \cup \Delta \sqsupset \Theta_2\} \Vdash \Delta_1 \cup \Delta_2 \sqsupset \Theta_1 \cup \Theta_2$.
\end{enumerate}
\end{lemma}

As elaborated before, over the class $\mathcal{D}_2$, the deductive system on i-formulas needs to be strengthened.
\begin{definition}[$\,\Vdash_2$]\label{definition of Vdash_2}
     For each formula $\gamma$, let $\,\Vdash^{\gamma}_2\, \subseteq \wp(\form^i) \times \form^i$ be the least $\,\Vdash$ that satisfies the basic i-rules, as well as the following additional rule:
    \begin{enumerate}
        \item[] $(\gamma$-Refl) \quad $\Vdash \{\varphi_1, \dots, \varphi_n, \bigvee_j(\varphi_j \rightarrow \psi_j) \lor \gamma\} \sqsupset \{\psi_1, \dots, \psi_n\}$.
    \end{enumerate}
\end{definition}

It is easily verified that $\Vdash^{\gamma}_2$ is compact, i.e., for each $\Gamma^{i}_{} \subseteq Form^{i}$ and $\varphi^{i} \in Form^{i}$, if $\Gamma^{i} \Vdash^{\gamma}_2 \varphi^{i}$, then there is a finite $\Delta^{i} \subseteq \Gamma^{i}$ such that $\Delta^{i} \Vdash^{\gamma}_2 \varphi^{i}$. The same holds for $\Vdash_1$ as well.

Next, we define a function that maps each set of formulas to a set of i-formulas: 
\begin{definition}
  For each $\gamma\in Form$ and $\Gamma\subseteq Form$, let
    $$i_{\gamma}(\Gamma)=\{\{\varphi_1, \dots, \varphi_n\} \sqsupset \{\psi_1, \dots, \psi_n\} \mid \bigvee\nolimits_j (\varphi_j \rightarrow \psi_j) \lor \gamma\in\Gamma\}.$$
\end{definition}

Now, we are at the point to define the syntactic consequence $\,\vdash_{k}$.
\begin{definition}[\(\, \vdash_{k} \)]
  Let \(\, \vdash_{1} \)  be the smallest relation $\,\vdash$ that satisfies all the basic rules (cf. Definition \ref{basic rule}), (Abs), ($\neg\neg$I) (cf. Definition \ref{def: (Abs) and (DNI)}) and the following rule:
  \begin{enumerate}
    \item[]
       $(Prop_{1})$ 
       \quad $i_{\gamma}(Th_{\,\vdash}(\alpha)) \Vdash_1 \{\top\} \sqsupset \{\bot\} \ \Rightarrow\;$ $\alpha \vdash \gamma$,
  \end{enumerate}
  where $Th_{\,\vdash}(\alpha)=\{\varphi\in Form \mid \alpha \vdash \varphi\}$. 
  Similarly,  let \(\, \vdash_{2} \)  be the smallest relation $\ \vdash$ that satisfies all the basic rules, ($\neg\neg$I), (Refl${}_1$), (Refl${}_2$) (cf. Definition \ref{def: (Refl1) and (Refl2)}) and the following rule:
  \begin{enumerate}
    \item[] 
       $(Prop_{2})$ 
       \quad $i_{\gamma}(Th_{\,\vdash}(\alpha)) \Vdash^{\gamma}_2 \{\top\} \sqsupset \{\bot\} \ \Rightarrow\;$ $\alpha \vdash \gamma$.
  \end{enumerate}
\end{definition}



The next lemma shows that for $\vdash_1$, the rule $(Prop_{1})$ is indeed necessary.
\begin{lemma}\label{necessity of Prop_1}
    $(Prop_{1})$ is not derivable from the other rules in the definition of $\,\vdash_1$.
\end{lemma}
\begin{proof}
    Let \(\, \vdash^{-}_{1} \) be the smallest relation  that satisfies all the defining rules of $\,\vdash_1$ except $(Prop_{1})$. From Lemma \ref{correctness part.1.1} and Lemma \ref{correspondance of (Abs) and (DNI)}, it readily follows that \( \, \vdash^{-}_{1} \) is sound with respect to the class of strongly pseudo-reflexive, weakly pseudo-symmetric $\Box\blackdiamond$-models. Now, let $p,q$ be two propositional letters, and let $\alpha$ be $((\top\to p)\lor q)\land ((p\to\bot)\lor q)$. Clearly, $i_{q}(Th_{\,\vdash^{-}_{1}}(\alpha)) \Vdash_1 \{\top\} \sqsupset \{\bot\}$. Below, we show $\alpha \nvdash^{-}_{1} q$ by giving a countermodel which is a strongly pseudo-reflexive, weakly pseudo-symmetric $\Box\blackdiamond$-model. It follows that \(\, \vdash^{-}_{1} \) fails $(Prop_{1})$.
    \begin{figure}[!ht]
        \centering
        \[
    \xymatrix@C=1pc@R=2.7pc{
     e: \textcolor{blue}{q} \ar@(ul,ur)[]  &  & &  & b: \textcolor{blue}{p}\ar@(ul,ur)[]\ar[dll]  \\ 
       & &  d: \textcolor{blue}{all} & &   \\
        a: \textcolor{blue}{none} \ar[uu]\ar[urr]  &  & &  & c: \textcolor{blue}{none}\ar@(dr,dl)[]\ar[ull]  }
    \]
  \end{figure}

    Consider the model $\M$ with 5 points as illustrated in the figure above, where the arrows represent the accessibility relation and the labels colored blue on the points indicate which propositional letters are true at that point (``all'' means all propositional letters are true there while ``none'' means none is true there).
    It is not hard to check that $\M$ is a strongly pseudo-reflexive, weakly pseudo-symmetric $\Box\blackdiamond$-model. Also, one can verify that $\M ,b \vDash \top\to p$ and $\M ,c \vDash p\to\bot$, which, together with $\M ,e \vDash q$, imply that
    $\M ,a \vDash ((\top\to p)\lor q)\land ((p\to\bot)\lor q)$. However, $\M ,a \nvDash q$. Thus, $\M$ is the desired countermodel.
\end{proof}

\section{Soundness and Completeness}
\subsection{Soundness}

\begin{definition}[Interpretation of i-formulas]
Let $\M = \langle \F, V \rangle $ be a Kripke model. 
We define the truth set of the i-formula $\Delta \sqsupset \Theta$ in $\M$:
\[
\Vert \Delta \sqsupset \Theta \Vert_{\M} = -_{\mathfrak{F}}\Vert \bigwedge \Delta \Vert_\M \cup \Vert \bigvee \Theta \Vert_\M .
\]
Additionally, for each $\Gamma^i \subseteq Form^i$, we write $\Vert \bigwedge \Gamma^i \Vert_{\M}$ for the intersection $\bigcap \{\Vert \alpha^i \Vert_{\M} \mid \alpha^i \in \Gamma^i\}$.
\end{definition}

It is not hard to verify that:
\begin{lemma}\label{soundness of vdash^i}
For each $\Gamma^i \subseteq Form^i$, $\varphi^i \in Form^i$, and $\gamma \in Form$, 
\begin{enumerate}
    \item if $\Gamma^i \Vdash_1 \varphi^i$, then \( \Vert \bigwedge \Gamma^i \Vert_{\M} \subseteq \Vert \varphi^i \Vert_{\M}\) for every Kripke model $\mathfrak{M}$;
    \item if $\Gamma^i \Vdash^{\gamma}_2 \varphi^i$, then $\Vert \bigwedge \Gamma^i \Vert_{\M} \subseteq \Vert \varphi^i \Vert_{\M} \cup \blackdiamond\Vert \gamma \Vert_{\M}$ for every reflexive Kripke model $\mathfrak{M}$.
\end{enumerate}

\end{lemma}

Observe that, for every Kripke model $\M$, $\|\bigvee(\varphi_j\to\psi_j)\lor\gamma\| = \Box_{\F}\blackdiamond_{\F}(\bigcup\|\varphi_j\to\psi_j\| \cup \|\gamma\|)$ $\subseteq$ $\Box_{\F} (\bigcup(-\|\varphi_j\| \cup \|\psi_j \|) \cup \blackdiamond_{\F}\|\gamma\|)$ $=\Box_{\F}(\|\{\varphi_1, \dots, \varphi_n\} \sqsupset \{\psi_1, \dots, \psi_n\}\| \cup \blackdiamond\|\gamma\|)$. Therefore, we have:

\begin{lemma}\label{semantic property of i_gamma(Gamma)}
 For each Kripke model $\M = \langle \F, V \rangle$, $\Gamma \subseteq Form$, and $\gamma \in Form$, 
    \[
    \Vert \bigwedge \Gamma \Vert_{\M} \subseteq \Box_{\F}(\Vert \bigwedge i_{\gamma}(\Gamma) \Vert_{\M} \cup \blackdiamond_{\F}\Vert \gamma \Vert_{\M}).
    \]
    ($\Vert \bigwedge \Gamma \Vert_{\M}$ denotes the intersection $\bigcap \{\Vert \alpha \Vert_{\M} \mid \alpha \in \Gamma\}$ when $\Gamma$ is infinite.)
\end{lemma}

\begin{lemma}\label{correctness of (Prop_1)}
    $(\text{Prop}_{1})$ is correct at every pseudo-symmetric $\Box\blackdiamond$-model. $(\text{Prop}_{2})$ is correct at every reflexive $\Box\blackdiamond$-model.
\end{lemma}
\begin{proof}
    We only prove the correctness of $(\text{Prop}_{2})$, as the proof of $(\text{Prop}_{1})$ is similar.
    Let $\M=<\F, V>$ be a reflexive $\Box\blackdiamond$-model. Assume that $i_{\gamma}(Th_{\,\vDash_{\M}}(\alpha)) \Vdash^{\gamma}_2 \{\top\} \sqsupset \{\bot\}$.
    Then by Lemma \ref{soundness of vdash^i}, we have $\Vert \bigwedge i_{\gamma}(Th_{\,\vDash_{\M}}(\alpha)) \Vert \subseteq \Vert \{\top\} \sqsupset \{\bot\} \Vert\cup \blackdiamond_{\F}\Vert \gamma \Vert = \Vert\bot\Vert\cup \blackdiamond_{\F}\Vert \gamma \Vert = \Box_{\F}\emptyset\cup \blackdiamond_{\F}\Vert \gamma \Vert$. Combined with Lemma \ref{semantic property of i_gamma(Gamma)}, we obtain  $\Vert \bigwedge Th_{\,\vDash_{\M}}(\alpha) \Vert \subseteq \Box_{\F}(\Box_{\F}\emptyset \cup \blackdiamond_{\F}\Vert \gamma \Vert)$. Then, since $\M$ is successor-serial (by reflexivity), and since $\Vert \gamma \Vert$ is a $\Box_{\F}\blackdiamond_{\F}$-fixpoint (by Lemma \ref{formula.fixedpoint}), we have $\Vert \bigwedge Th_{\,\vDash_{\M}}(\alpha) \Vert$ $\subseteq \Box_{\F}\blackdiamond_{\F}\Vert \gamma \Vert \subseteq\Vert \gamma \Vert$. Clearly, $\Vert\alpha\Vert \subseteq \bigcap \{\Vert \varphi \Vert \mid \alpha \vDash_{\M}\varphi\} = \Vert \bigwedge Th_{\,\vDash_{\M}}(\alpha) \Vert$, thus, it follows that $\alpha \vDash_{\M} \gamma$.
\end{proof}

\begin{theorem}
For every $\Gamma\subseteq Form$, $\varphi\in Form$, and $\ k\in\{1,2\}$,
    \[\Gamma\,\vdash_{k}\varphi \Rightarrow\Gamma\,\vDash_{k}\varphi.\]
\end{theorem}

\begin{proof}
    We only prove the case of \(\vdash_{2}\), as the case of \( \vdash_{1}\) is similar.
    Let $\mathfrak{M}\in \mathcal{D}_{2}$. In order to prove $\,\vdash_{2}\, \subseteq\, \,\vDash_{\mathfrak{M}}$, it suffices to show that $\,\vDash_{\mathfrak{M}}$ satisfies all the basic rules, ($\neg\neg$I), (Refl${}_1$), (Refl${}_2$), and (Prop${}_{2}$)---i.e., all the defining rules of $\,\vdash_{2}$ are correct at $\M$. This follows from  Lemma \ref{correctness part.1.1}, Corrollary \ref{correctness part.1.2}, Lemma \ref{correctness part.1.3}, and Lemma \ref{correctness of (Prop_1)}.
    \qed
\end{proof}

\subsection{Completeness}\label{sec: completeness}
To prove the completeness of $\vdash_{k}$, our approach is to build a canonical model $\M^c_k$ using pairs of sets of formulas such that $\M^c_k\in\mathcal{D}_k$, and for each  $\langle \Gamma, \Delta \rangle\in\M^c_k$, $\Gamma$ is the set of formulas that are true at $\langle \Gamma, \Delta \rangle$, while $\Delta$ collects those formulas that are  false at every predecessor of $\langle \Gamma, \Delta \rangle$. And the accessibility relation is defined by: $\langle \Gamma_1, \Delta_1 \rangle R^c_k \langle \Gamma_2, \Delta_2 \rangle \Leftrightarrow \Gamma_1 \cap \Delta_2 = \emptyset$.

Towards this end, for each point $\langle \Gamma, \Delta \rangle$ in $\M^c_k$, the following conditions should be satisfied:
\begin{enumerate}[(1)]
    \item 
    $\Delta$ contains $\bot$ and is closed under disjunction;
    \item for any $\alpha \in \text{Form}$, $(\Gamma \vDash_{k} \alpha \Rightarrow \alpha \in \Gamma)$; 
    \item 
    for any $\alpha,\beta \in \text{Form}$, $(\alpha \vDash_{k} \beta \text{ and } \beta \in \Delta \ \Rightarrow\ \alpha \in \Delta)$;
    \item $\{\alpha \rightarrow \beta \mid \alpha \in \Gamma$ and $\beta \notin \Gamma\}\subseteq\Delta$;
    \item $\{\neg\alpha \mid \alpha \in \Gamma\}\subseteq\Delta$;
    \item only for $\M^c_2$: $\Gamma\cap\Delta=\emptyset$.
\end{enumerate}
(For example, given that $\Delta$ is the set of formulas false at every predecessor, condition (1) and (4) follows from the semantic definitions of $\bot$ and $\lor$, while condition (5) follows from Lemma \ref{semantics: neg = to bot} and the successor-seriality of $\M^c_k$.)

\begin{definition}
Given $\,\vdash\,\subseteq \wp(Form)\times Form$, for each $\Gamma,\Delta$, we define:
   \begin{itemize}
       \item  $\Gamma$ is $\,\vdash$-closed  iff $Th_{\,\vdash}(\Gamma)\subseteq\Gamma$, where $Th_{\,\vdash}(\Gamma)=\{\alpha\in Form\mid \Gamma\vdash\alpha\}$.  
       \item  $\Delta$ is $\,\vdash$-downward closed  iff for each $\alpha,\beta$, $(\beta\in\Delta$ and $\alpha\vdash \beta \Rightarrow \alpha \in \Delta)$.
       
       \item $\neg(\Gamma)=\{\neg\alpha \mid \alpha \in \Gamma\}$, and $\!\to\!(\Gamma)=\{\alpha \rightarrow \beta \mid \alpha \in \Gamma,\ \beta \notin \Gamma\}$. 
   \end{itemize}
\end{definition}

As elaborated above,  the canonical model for $\vdash_k$ should be a submodel of the model $W^d_{k}$ defined below.
\begin{definition}
    Define \(\mathfrak{M}^d_1 = \langle W^d_1, R^d_1, V^d_1 \rangle\) as follows:
    \begin{itemize}
        \item $W^d_1=\{\langle \Gamma, \Delta \rangle\in (\wp(Form))^2  \mid \Gamma$ is $\,\vdash_1$-closed,  $\{\bot\}\cup\neg(\Gamma)\cup\!\to\!(\Gamma)\subseteq\Delta$, and $\Delta$ is $\,\vdash_1$-downward closed and closed under disjunction$\}.$

        \item \(R^d_1 = \{\langle \langle \Gamma_1, \Delta_1 \rangle, \langle \Gamma_2, \Delta_2 \rangle \rangle \in (W^d_1)^2 \mid \Gamma_1 \cap \Delta_2 = \emptyset\}\);
    
        \item \(V^d_1\) maps each $p\in PL$ to \( \{\langle \Gamma, \Delta \rangle \in W^d_1 \mid p \in \Gamma\}\).
    \end{itemize}
     \(\mathfrak{M}^d_2 = \langle W^d_2, R^d_2, V^d_2 \rangle\) is defined similarly, except that:
    \begin{itemize}
        \item $W^d_2=\{\langle \Gamma, \Delta \rangle\in (\wp(Form))^2  \mid \Gamma$ is $\,\vdash_2$-closed,  $\{\bot\}\cup\neg(\Gamma)\cup\!\to\!(\Gamma)\subseteq\Delta$, $\Delta$ is $\,\vdash_2$-downward closed and closed under disjunction, and $\Gamma\cap\Delta=\emptyset\}.$
    \end{itemize}   
\end{definition}

Next, let us look at a specific kind of element within $W^d_k$.

\begin{definition}\label{definition of I_gamma(Gamma)}
For each $k\in\{1,2\}$, $\gamma\in Form$ and each $\Gamma\subseteq Form$, 
let
$I_k^{\gamma}(\Gamma) = \{\varphi \mid$ there exists $\alpha_1,\dots,\alpha_n \in \Gamma$ and $\beta_1, \dots, \beta_n \notin  \Gamma$ (for some $n\in\mathbb{N}^+$) such that $\varphi \vdash_k (\alpha_1 \rightarrow \beta_1) \lor \dots \lor (\alpha_n \rightarrow \beta_n) \lor \gamma\}$.
Then define 
\begin{itemize}
    \item $W^b_{1}=\{\langle \Gamma, I_1^{\gamma}(\Gamma) \rangle \mid  \Gamma\text{ is } \vdash_1 \text{-closed and }\bot\notin\Gamma\}$;
    \item $W^b_{2}=\{\langle \Gamma, I_2^{\gamma}(\Gamma) \rangle \mid  \Gamma\text{ is } \vdash_2 \text{-closed, }\bot\notin\Gamma, \text{ and } \Gamma\cap I_2^{\gamma}(\Gamma)=\emptyset\}$.
\end{itemize}
\end{definition} 

Using the $\to$-free basic rules only, one can show the following lemma. (Note that the condition ``$\bot\notin\Gamma$'' in the defintion of $W^b_{k}$ is used when proving $\neg(\Gamma)\subseteq I_k^{\gamma}(\Gamma)$.)
\begin{lemma}
    $W^b_{k}\subseteq W^d_{k}$, for each $k\in\{1,2\}$.
\end{lemma}

Next, we are going to prove five existence lemmas: Extension Lemma, Implication Lemma, Fixpoint Lemma, Pseudo-reflexivity Lemma, and Pseudo-symmetry Lemma. They are crucial for establishing the completeness theorem later. But first, let us introduce the following encompassing lemma that will significantly simplify the proofs of these existence lemmas.


\begin{lemma}\label{the encompassing lemma}
    Let  $\,\vdash$ and $\,\vsim$ be binary relations on formulas such that $\,\vdash\,\subseteq\, \vsim$ and $\,\vsim$ satisfies (A), (Cut), ($\land$L/R) and ($\lor$R). 
    Let $F$ map each finite sequence of pairs $\langle\alpha_j,\beta_j\rangle$ 
    to a formula. Define $I_{\,\vdash,F}(\Delta)=\{\varphi \mid \text{there exists } \alpha_1,\dots,\alpha_n \in \Delta$ and $\beta_1, \dots, \beta_n \notin  \Delta$ with $\varphi \vdash\! F(\langle\alpha_1 ,\beta_1\rangle , \dots , \langle\alpha_n ,\beta_n\rangle)\}$. Let $I\subseteq Form$ be 
    either a singleton or closed under disjunction, and let $\mathcal{S}_{I}=\{\Delta\subseteq Form \mid \Delta$ is $\vsim$-closed and $\Delta\cap I=\emptyset\}$. If $\mathcal{S}$ is a end-segment of $\mathcal{S}_{I}$ (i.e., $\mathcal{S}\subseteq \mathcal{S}_{I}$ and for each $\Delta_1,\Delta_2$, if $\Delta_1\subseteq\Delta_2$ and $\Delta_1\in \mathcal{S}$, then $\Delta_2\in \mathcal{S}$), and if $\Phi$ is maximal in $\mathcal{S}$, then:
    \begin{enumerate}
        \item For each $\Gamma\subseteq Form$, if $\,\vsim$ satisfies the following rule:\\
        ($\Gamma\vdash\! F$)\quad   $\Gamma\vdash\! F(\langle\alpha_1 ,\beta_1\rangle , \dots , \langle\alpha_n ,\beta_n\rangle)$ and $\psi_j \land \beta_j \vsim \varphi$ for $1\leq j\leq n$ \\
        \indent\qquad\qquad~$\Rightarrow \bigwedge_j(\psi_j \land \alpha_j) \vsim \varphi$\\
        then $\Gamma \cap I_{\,\vdash,F}(\Phi) = \emptyset$.

        \item If $\,\vsim$ satisfies the following rule:\\
             ($F$-MP)\quad  $\psi_j \land \beta_j \vsim \varphi$ for $1\leq j\leq n$ 
             $\Rightarrow \bigwedge_j(\psi_j \land \alpha_j) \land F(\langle\alpha_1 ,\beta_1\rangle , \dots , \langle\alpha_n ,\beta_n\rangle)$\\
             \indent\qquad\qquad~$\vsim \varphi$\\
        then $\Phi \cap I_{\,\vdash,F}(\Phi) = \emptyset$.

    \end{enumerate}
\end{lemma}

\begin{proof}
    \begin{enumerate}
        \item Let $\Gamma\subseteq Form$. Assume $\,\vsim$ satisfies ($\Gamma\vdash\! F$). Suppose, towards a contradiction, that $\Gamma \cap I_{\,\vdash,F}(\Phi)$ is non-empty. Then, there exist $\alpha_1,\dots,\alpha_n \in \Phi$ and $\beta_1, \dots, \beta_n \notin  \Phi$ such that $\Gamma \vdash\! F(\langle\alpha_1 ,\beta_1\rangle , \dots , \langle\alpha_n ,\beta_n\rangle)$. for each $1\leq j\leq n$, since $\beta_j \notin \Phi$ and $\Phi$ is maximal in $\mathcal{S}$, we have $Th_{\vsim}(\Phi \cup \{\beta_j\}) \notin \mathcal{S}$(---since $\vsim$ satisfies (A), $\Phi \cup \{\beta_j\}\subseteq Th_{\vsim}(\Phi \cup \{\beta_j\})$); then, since $\mathcal{S}$ is a end-segment of $\mathcal{S}_{I}$, and since $Th_{\vsim}(\Phi \cup \{\beta_j\})$ is $\vsim$-closed(---using the fact that $\vsim$ satisfies (Cut) and ($\land$L/R)), we have $Th_{\vsim}(\Phi \cup \{\beta_j\})\cap I\neq\emptyset$. Then, we can find a uniform $\chi\in I$ with $\chi\in \bigcap_j Th_{\vsim}(\Phi \cup \{\beta_j\})$: this is obvious when $I$ is a singleton; otherwise, if $I$ is closed under disjunction, suppose for each $j$, $\chi_j\in Th_{\vsim}(\Phi \cup \{\beta_j\})\cap I$, then $\chi:=\bigvee_j \chi_j$ belongs to $I\cap\bigcap_j Th_{\vsim}(\Phi \cup \{\beta_j\})$(---using the fact that $\vsim$ satisfies ($\lor$R) and (Cut)). 
        
        It follows, from $\chi\in\bigcap_j Th_{\vsim}(\Phi \cup \{\beta_j\})$, that for each $j$, there exists finite $\Phi_j \subseteq \Phi$ such that $\bigwedge\Phi_j \land \beta_j \vsim \chi$(---using the fact that $\vsim$ satisfies (Cut) and ($\land$L/R)).
        By ($\Gamma\vdash\! F$), we have $\bigwedge_j(\bigwedge\Phi_j \land \alpha_j) \vsim \chi$. Since $\Phi$ is $\vsim$-closed and $\bigcup_j(\Phi_j\cup\{\alpha_j\})$ is a finite subset of $\Phi$, it follows that $\chi\in\Phi$, thus, $\Phi \cap I\neq\emptyset$, contradicting the fact that $\Phi\in\mathcal{S}_{I}$.
        
        \item Assume $\,\vsim$ satisfies ($F$-MP). Suppose, towards a contradiction, that $\Phi \cap I_{\,\vdash,F}(\Phi)$ is non-empty. Then, there exist $\alpha_1,\dots,\alpha_n \in \Phi$ and $\beta_1, \dots, \beta_n \notin  \Phi$ such that $\Phi \vdash\! F(\langle\alpha_1 ,\beta_1\rangle , \dots , \langle\alpha_n ,\beta_n\rangle)$, thus, by $\,\vdash\,\subseteq\, \vsim$ we have  $\Phi \vsim F(\langle\alpha_1 ,\beta_1\rangle , \dots , \langle\alpha_n ,\beta_n\rangle)$. Moreover, since $\beta_j \notin \Phi$ for each $j$, similarly to the proof of (1),
        we can find $\chi\in I$ and finite subsets $\Phi_1,\dots, \Phi_n$ of $\Phi$ such that $\bigwedge\Phi_j \land \beta_j \vsim \chi$ for each $ j$. By ($F$-MP), we have $\bigwedge_j(\bigwedge\Phi_j \land \alpha_j) \land F(\langle\alpha_1 ,\beta_1\rangle , \dots , \langle\alpha_n ,\beta_n\rangle)\vsim \chi$. Since $\Phi$ is $\vsim$-closed, 
        it follows that $\chi\in\Phi$, thus, $\Phi \cap I\neq\emptyset$, contradicting the fact that $\Phi\in\mathcal{S}_{I}$.

    \end{enumerate}
\end{proof}

In order to apply the above lemma to prove the five existence lemmas, we consider the following three instances of $\vsim$:
\begin{lemma}\label{three instances of vsim}
    Let $\,\vdash$ be a consequence relation on formulas satisfying the basic rules, and let $\,\Vdash$ be a consequence relation on i-formulas satisfying the basic i-rules (cf. Definition \ref{definition of Vdash_1}). Fix $\Gamma\subseteq Form$ and $\gamma\in Form$, define the following relations:
    \begin{itemize}
        \item  $\alpha \vsim_1 \beta$ $\iff$ $\alpha \vdash   \beta$;
        \item  $\alpha \vsim_2 \beta$ $\iff$ $\Gamma \vdash \alpha \to \beta$;
        \item  $\alpha \vsim_3 \beta$ $\iff$ $i_{\gamma}(Th_{\,\vdash}(\Gamma)) \Vdash \{\alpha\} \sqsupset \{\beta\}$.        
    \end{itemize}
    Then, for each $j\in\{1,2,3\}$, $\,\vdash\,\subseteq\, \vsim_j$, and $\vsim_j$ satisfies (A), (Cut), ($\bot$), ($\land$L/R) and ($\lor$R). 
    Moreover, define $f_{\gamma}(\langle\alpha_1 ,\beta_1\rangle , \dots , \langle\alpha_n ,\beta_n\rangle) = \bigvee_j(\alpha_j \rightarrow \beta_j) \lor \gamma$ for each finite sequence of pairs $\langle\alpha_j,\beta_j\rangle$.
    Then:
    \begin{enumerate}
        \item If $\,\vdash$ satisfies (Refl${}_1$) (cf. Definition \ref{def: (Refl1) and (Refl2)}), then $\vsim_1$ satisfies ($f_{\bot}$-MP).
        
        \item 
          \begin{enumerate}
            \item $\vsim_2$ satisfies ($\Gamma\vdash\! f_{\bot}$).
            \item If $\,\vdash$ satisfies (Refl${}_2$), then $\,\vsim_2$ satisfies ($f_{\bot}$-MP).
          \end{enumerate}  
        
        \item 
          \begin{enumerate}
            \item $\vsim_3$ satisfies ($\Gamma\vdash\! f_{\gamma}$).
            \item If $\,\Vdash$ satisfies ($\gamma$-Refl) (cf. Definition \ref{definition of Vdash_2}), then $\,\vsim_3$ satisfies ($f_{\gamma}$-MP).
          \end{enumerate}  
    \end{enumerate}    
\end{lemma}

\begin{proof}
    It is straightforward to verify that for each $j\in\{1,2,3\}$, $\,\vdash\,\subseteq\, \vsim_j$ and $\vsim_j$ satisfies (A), (Cut), ($\bot$), ($\land$L/R), and ($\lor$R). Note that $\vsim_j$ satisfying (A), ($\bot$), ($\land$L) and ($\lor$R) is implied by $\,\vdash$ satisfying them together with $\,\vdash\,\subseteq\, \vsim_j$. 
    Below, we prove the clauses 1$\sim$3 respectively.
    \begin{enumerate}
        \item Assume $\,\vdash$ satisfies (Refl${}_1$). Suppose $\psi_j \land \beta_j \vsim_1 \varphi$, i.e., $\psi_j \land \beta_j \vdash \varphi$, for $1\leq j\leq n$. By (Refl${}_1$), we have $\bigwedge_j(\psi_j \land \alpha_j) \land \bigvee_j(\alpha_ \rightarrow \beta_j) \vdash \varphi$. Thus, $\bigwedge_j(\psi_j \land \alpha_j) \land f_{\bot}(\langle\alpha_1 ,\beta_1\rangle , \dots , \langle\alpha_n ,\beta_n\rangle)\vsim_1 \varphi$.
        
        \item 
          \begin{enumerate}
            \item This follows from the fact that $\,\vdash$ satisfies ($\to$-$\lor$.s), ($\bot$), ($\lor$L) etc.
            \item Assume $\,\vdash$ satisfies (Refl${}_2$). Suppose $\psi_j \land \beta_j \vsim_2 \varphi$, i.e., $\Gamma \vdash \psi_j \land \beta_j \to \varphi$, for $1\leq j\leq n$. By (Refl${}_2$), ($\land$L/R) and (Cut), we have $\Gamma \vdash \bigwedge_j(\psi_j \land \alpha_j) \land \bigvee_j(\alpha_j \rightarrow \beta_j) \to \varphi$. Thus, $\bigwedge_j(\psi_j \land \alpha_j) \land f_{\bot}(\langle\alpha_1 ,\beta_1\rangle , \dots , \langle\alpha_n ,\beta_n\rangle)\vsim_2 \varphi$.
          \end{enumerate}

        \item 
          \begin{enumerate}
            \item If $\Gamma \vdash f_{\gamma}(\langle\alpha_1 ,\beta_1\rangle , \dots , \langle\alpha_n ,\beta_n\rangle)$ and $\psi_j \land \beta_j \vsim_3 \varphi$ for $1\leq j\leq n$, then by the definition of $i_{\gamma}$ and $\vsim_3$, we have $\{\alpha_1,\dots,\alpha_n\}\sqsupset \{\beta_1,\dots,\beta_n\}\in i_{\gamma}(Th_{\,\vdash}(\Gamma))$, and $i_{\gamma}(Th_{\,\vdash}(\Gamma)) \Vdash \{\psi_j \land \beta_j\} \sqsupset \{\varphi\}$ for $1\leq j\leq n$. 
            From the former, using basic i-rules we can prove that $i_{\gamma}(Th_{\,\vdash}(\Gamma)) \Vdash\{\bigwedge_j(\psi_j \land \alpha_j) \} \sqsupset \{\psi_1 \land \beta_1,\, \dots\,,\, \psi_n \land \beta_n \}$; together with the latter, using (i-Cut) we get $i_{\gamma}(Th_{\,\vdash}(\Gamma)) \Vdash \{\bigwedge_j(\psi_j \land \alpha_j) \} \sqsupset \{\varphi\}$. Thus, $\bigwedge_j(\psi_j \land \alpha_j) \vsim_3 \varphi$.
            
            \item Assume $\,\Vdash$ satisfies ($\gamma$-Refl). Suppose $\psi_j \land \beta_j \vsim_3 \varphi$, i.e., $i_{\gamma}(Th_{\,\vdash}(\Gamma)) \Vdash \{\psi_j \land \beta_j\} \sqsupset \{\varphi\}$, for $1\leq j\leq n$. By ($\gamma$-Refl) and basic i-rules, we can get $\Vdash \{\bigwedge_j(\psi_j \land \alpha_j) \land f_{\gamma}(\langle\alpha_1 ,\beta_1\rangle , \dots , \langle\alpha_n ,\beta_n\rangle) \} \sqsupset \{\psi_1 \land \beta_1,\, \dots\,,\, \psi_n \land \beta_n\}$. Together with $i_{\gamma}(Th_{\,\vdash}(\Gamma)) \Vdash \{\psi_j \land \beta_j\} \sqsupset \{\varphi\}$, we have 
            $i_{\gamma}(Th_{\,\vdash}(\Gamma)) \Vdash \{\bigwedge_j(\psi_j \land \alpha_j) \land f_{\gamma}(\langle\alpha_1 ,\beta_1\rangle , \dots , \langle\alpha_n ,\beta_n\rangle) \} \sqsupset \{\varphi\}.$
            Thus, $\bigwedge_j(\psi_j \land \alpha_j) \land f_{\gamma}(\langle\alpha_1 ,\beta_1\rangle , \dots , \langle\alpha_n ,\beta_n\rangle)\vsim_3 \varphi$.
          \end{enumerate}
    \end{enumerate}
\end{proof}

Now, we are ready to prove the five existence lemmas mentioned before.
\begin{lemma}\label{five existence lemmas}
    Let $k\in\{1,2\}$.
\begin{enumerate}
    \item \textbf{Extension Lemma}\\
          If $ I\subseteq Form$ is nonempty and closed under disjunction, and if $Th_{\,\vdash_k}(\Gamma)\cap I=\emptyset$, then there exists $\langle \Phi, I^{\bot}_k(\Phi) \rangle\in W^b_{k}$ such that $\Gamma \subseteq \Phi \text{ and }  \Phi\cap I=\emptyset$.   
    \item \textbf{Implication Lemma}\\
          If $\Gamma\nvdash_k\alpha \to \beta$, then there is $\langle \Phi, I^{\bot}_k(\Phi) \rangle\in W^b_{k}$ such that $\Gamma \cap I^{\bot}_k(\Phi) = \emptyset \text{, } a \in \Phi \text{ and } \beta \notin \Phi$.
    \item \textbf{Fixpoint Lemma}\\
          If $\Gamma\nvdash_k\gamma $, then there exists  $\langle \Phi, I_{k}^{\gamma}(\Phi) \rangle\in W^b_{k}$ such that $\Gamma \cap I_{k}^{\gamma}(\Phi) = \emptyset$.
    \item \textbf{Pseudo-reflexivity Lemma}\\
          For each $\langle\Gamma,\Delta\rangle \in W^d_{1}$, if $\bot\notin\Gamma$, then there exists \(\langle \Phi, I^{\bot}_1(\Phi) \rangle\in W^b_{1}\) such that \(\Gamma \subseteq \Phi \text{ and } \Gamma \cap I^{\bot}_1(\Phi) = \emptyset.\)
    \item \textbf{Pseudo-symmetry Lemma}\\
          For each $\langle\Gamma,\Delta\rangle, \langle\Gamma',\Delta'\rangle \in W^d_{k}$, if $\Gamma \cap \Delta' = \emptyset$, then \( \text{there is } \langle \Phi, I^{\bot}_k(\Phi) \rangle\in W^b_{k}\) such that \(\Gamma \subseteq \Phi \text{ and } \Gamma' \cap I^{\bot}_k(\Phi) = \emptyset.\)
\end{enumerate}
\end{lemma}




\begin{proof}
    \begin{enumerate}
        \item Assume $Th_{\,\vdash_k}(\Gamma)\cap I=\emptyset$. For $k=1$, we are done by simply taking $\Phi=Th_{\,\vdash_1}(\Gamma)$---note that by ($\bot$) and (Cut), it follows from $I\neq\emptyset$ and $Th_{\,\vdash_1}(\Gamma)\cap I=\emptyset$ that $\bot\notin Th_{\,\vdash_1}(\Gamma)$. 
        For $k=2$, by Lemma \ref{three instances of vsim}, we know that $\vdash_2$ satisfies ($f_{\bot}$-MP).
        Let $\mathcal{S}=\{\Psi\subseteq Form \mid \Psi$ is $\vdash_2$-closed, $\Gamma\subseteq\Psi$,  $\Psi\cap I=\emptyset\}$. 
        Clearly, $Th_{\vdash_2}(\Gamma) \in \mathcal{S}$. Moreover, for any chain $H \subseteq \mathcal{S}$, it is easy to see that $\bigcup H \in \mathcal{S}$. 
        So, by Zorn's lemma, $\mathcal{S}$ has a maximal element, denoted by $\Phi$.  
        By  Lemma \ref{the encompassing lemma}, we have $\Phi \cap I^{\bot}_2(\Phi) = \emptyset$. Additionally, by ($\bot$) and (Cut), it follows from $I\neq\emptyset$ and $\Phi\cap I=\emptyset$ that $\bot\notin \Phi$. Thus, $\langle \Phi, I^{\bot}_2(\Phi) \rangle\in W^b_{2}$ is as required.
        
        \item Assume $\Gamma\nvdash_k\alpha \to \beta$. Define the relation $\vsim$: $\alpha\vsim\beta\iff\Gamma\vdash_k\alpha\to\beta$. By Lemma \ref{three instances of vsim}, we know that $\,\vdash\,\subseteq\, \vsim$, and $\vsim$ satisfies (A), (Cut), ($\bot$), ($\land$L/R), ($\lor$R), ($\Gamma\vdash_k\! f_{\bot}$), and (if $k=2$) ($f_{\bot}$-MP).
        
        Let $\mathcal{S}=\{\Psi\subseteq Form \mid \Psi$ is $\vsim$-closed, $\alpha \in \Psi$, $\beta\notin\Psi\}$. 
        Since $\vsim$ satisfies (A), (Cut) and ($\land$L/R), and since $\Gamma\nvdash_k\alpha \to \beta$, we can prove that $Th_{\vsim}(\{\alpha\}) \in \mathcal{S}$. Moreover, for any chain $H \subseteq \mathcal{S}$, it is easy to see that $\bigcup H \in \mathcal{S}$. 
        So, by Zorn's lemma, $\mathcal{S}$ has a maximal element, denoted by $\Phi$.  
        By  Lemma \ref{the encompassing lemma}, we have $\Gamma \cap I^{\bot}_k(\Phi) = \emptyset$ and (if $k=2$) $\Phi \cap I^{\bot}_k(\Phi) = \emptyset$. Moreover, since $\vsim$ satisfies ($\bot$), and since $\Phi$ is $\vsim$-closed and $\beta\notin\Phi$, we get $\bot\notin \Phi$. Thus, $\langle \Phi, I^{\bot}_k(\Phi) \rangle\in W^b_{k}$ is as required.

        \item Assume  $\Gamma \nvdash_k \gamma$. For convenience, let us denote $\Vdash_1$ by $\Vdash^\gamma_1$.  Define the relation $\vsim$: $\alpha\vsim\beta\iff i_{\gamma}(Th_{\,\vdash}(\Gamma)) \Vdash^\gamma_k \{\alpha\} \sqsupset \{\beta\}$. By Lemma \ref{three instances of vsim}, we know that $\,\vdash\,\subseteq\, \vsim$, and $\vsim$ satisfies (A), (Cut), ($\bot$), ($\land$L/R), ($\lor$R), ($\Gamma\vdash_k\! f_{\gamma}$), and (if $k=2$) ($f_{\gamma}$-MP).
        Let $\mathcal{S}=\{\Psi\subseteq Form \mid \Psi$ is $\vsim$-closed, $\gamma\notin\Psi\}$. 

        First, we show that $\bot\notin Th_{\vsim}(\emptyset)$. Suppose not. Then $\top\vsim\bot$, i.e., $i_{\gamma}(Th_{\,\vdash}(\Gamma))$ $\Vdash^\gamma_k \{\top\} \sqsupset \{\bot\}$. By (Prop${}_k$),  we can prove that  $\Gamma \vdash_k  \gamma$. Contradiction.

        Moreover, it is easy to see that for any chain $H \subseteq \mathcal{S}$,  $\bigcup H \in \mathcal{S}$. Then, using Zorn's lemma we can get  a maximal element $\Phi$ of $\mathcal{S}$, which, by Lemma \ref{the encompassing lemma}, satisfies that $\Gamma \cap I^{\gamma}_k(\Phi) = \emptyset$ and (if $k=2$) $\Phi \cap I^{\gamma}_k(\Phi) = \emptyset$. Moreover, since $\vsim$ satisfies ($\bot$), and since $\Phi$ is $\vsim$-closed and $\gamma\notin\Phi$, we get $\bot\notin \Phi$. Thus, $\langle \Phi, I^{\gamma}_k(\Phi) \rangle\in W^b_{k}$ is as required.

        \item Assume $\langle\Gamma,\Delta\rangle \in W^d_{1}$ and $\bot\notin\Gamma$. Define the relation $\vsim$ similarly to the proof of 2, except that, here, $k=1$. 
        Let $\mathcal{S}=\{\Psi\subseteq Form \mid \Psi$ is $\vsim$-closed, $\Gamma \subseteq \Psi$, $\bot\notin\Psi\}$. 

        First, we show that $\bot\notin Th_{\vsim}(\Gamma)$. Suppose not. Then, there exists finite $\Theta\subseteq\Gamma$ s.t. $\bigwedge\Theta\vsim\bot$, i.e., $\Gamma\vdash_1\neg\bigwedge\Theta$. Since $\Gamma$ is $\,\vdash_1$-closed (as $\langle\Gamma,\Delta\rangle \in W^d_{1}$), we have $\bigwedge\Theta\in\Gamma$ and $\neg\bigwedge\Theta\in\Gamma$, then by (Abs), we have $\bot\in\Gamma$. Contradiction.

        Then, analogous to the proof of 2, we know that $\mathcal{S}$ has a maximal element $\Phi$ s.t. $\Gamma \cap I^{\bot}_1(\Phi) = \emptyset$. Thus, $\langle \Phi, I^{\bot}_1(\Phi) \rangle\in W^b_{1}$ is as required.

        \item Assume $\langle\Gamma,\Delta\rangle, \langle\Gamma',\Delta'\rangle \in W^d_{k}$ with $\Gamma \cap \Delta' = \emptyset$. Define the relation $\vsim$: $\alpha\vsim\beta\iff\Gamma'\vdash_k\alpha\to\beta$. By Lemma \ref{three instances of vsim}, we know that $\,\vdash\,\subseteq\, \vsim$, and $\vsim$ satisfies (A), (Cut), ($\bot$), ($\land$L/R), ($\lor$R), ($\Gamma'\vdash_k\! f_{\bot}$), and (if $k=2$) ($f_{\bot}$-MP).
        Let $\mathcal{S}=\{\Psi\subseteq Form \mid \Psi$ is $\vsim$-closed, $\Gamma \subseteq \Psi$, $\bot\notin\Psi\}$. 

        First, we show that $\bot\notin Th_{\vsim}(\Gamma)$. Suppose not. Then, there exists finite $\Theta\subseteq\Gamma$ s.t. $\bigwedge\Theta\vsim\bot$, i.e., $\Gamma'\vdash_k\neg\bigwedge\Theta$. On one hand, by ($\lnot \lnot$I) and since both $\Gamma$ and $\Gamma'$ are $\,\vdash_k$-closed, we have $\neg\neg\bigwedge\Theta\in\Gamma$ and $\neg\bigwedge\Theta\in\Gamma'$. On the other hand, since $\langle\Gamma',\Delta'\rangle \in W^d_{\,\vdash}$,  we have $\neg(\Gamma')\subseteq\Delta'$, then, by $\Gamma \cap \Delta' = \emptyset$ we have $\Gamma \cap \neg(\Gamma') = \emptyset$. Contradiction.

        Then, using Zorn's lemma we can get  a maximal element $\Phi$ of $\mathcal{S}$, which, by Lemma \ref{the encompassing lemma}, satisfies that $\Gamma' \cap I^{\bot}_k(\Phi) = \emptyset$ and (if $k=2$) $\Phi \cap I^{\bot}_k(\Phi) = \emptyset$. Thus, $\langle \Phi, I^{\bot}_k(\Phi) \rangle\in W^b_{k}$ is as required.
    \end{enumerate}
\end{proof}

\begin{lemma}[Main lemma for the canonical model]\label{Main lemma for the canonical model}
Let $k\in\{1,2\}$. Take any $W$ such that $W^b_{k} \subseteq W \subseteq W^d_{k}$, and let $\mathfrak{M} = \langle\mathfrak{F},V\rangle= \langle W,R,V\rangle$ be the submodel of $\mathfrak{M}^d_{k}$ on $W$. Define 
 $\vert \alpha\vert _W = \{\langle\Gamma, \Delta\rangle \in W \mid \alpha \in \Gamma\}$ for each formula $\alpha$. Then we have:
\begin{enumerate}[(1)]
    \item 
    For each $\alpha \in Form$, $\vert \alpha\vert _{W} \in FP_{\mathfrak{F}}$.
    \item \textbf{Truth Lemma:} For each $\alpha \in Form$, $\vert \alpha\vert _{W} = \Vert\alpha\Vert_{\mathfrak{M}}$.
    \item  $\mathfrak{F}$ is pseudo-symmetric. If $k=1$, then $\mathfrak{F}$ is pseudo-reflexive. If $k=2$, then $\mathfrak{F}$ is reflexive.
\end{enumerate}
\end{lemma}

\begin{proof}
(In this proof, we omit the subscript $\mathfrak{M}$ from $\Vert\cdot\Vert$, and the subscript $W$ from $\vert \cdot\vert $.)
\begin{enumerate}[(1)]
    \item Let $\alpha \in Form$. We show that $\Box_{\mathfrak{F}} \blackdiamond_{\mathfrak{F}}\vert \alpha\vert  \subseteq \vert \alpha\vert $, which is equivalent to $-_{\mathfrak{F}}\vert \alpha\vert  \subseteq \Diamond_{\mathfrak{F}}-_{\mathfrak{F}}\blackdiamond_{\mathfrak{F}}\vert \alpha\vert $. For any $\langle\Gamma, \Delta\rangle \in W$, if $\langle\Gamma, \Delta\rangle \in -_{\mathfrak{F}}\vert \alpha\vert $, i.e., $\alpha \notin \Gamma$, then, by Fixpoint Lemma (i.e., Lemma \ref{five existence lemmas}.3), there exists $\langle\Phi, I_{k}^{\alpha}(\Phi)\rangle \in W^b_{k}\subseteq W$ such that $\langle\Gamma, \Delta\rangle R \langle\Phi, I_{k}^{\alpha}(\Phi)\rangle$. Since $\alpha\in I_{k}^{\alpha}(\Phi)$, it follows that $\langle\Phi, I_{k}^{\alpha}(\Phi)\rangle \in -_{\mathfrak{F}}\blackdiamond_{\mathfrak{F}}\vert \alpha\vert $.
    
    \item 
    First, we prove that
    for any $\alpha, \beta \in Form$, we have:
    \begin{itemize}
        \item $\vert \bot\vert  = \Box_{\mathfrak{F}}\emptyset$ 
        \item $\vert \alpha \to \beta\vert  = \Box_{\mathfrak{F}}(-_{\mathfrak{F}}\vert \alpha\vert \cup \vert \beta\vert )$
        \item $\vert \alpha \land \beta\vert  = \vert \alpha\vert  \cap \vert \beta\vert $
        \item $\vert \alpha \lor \beta\vert  = \Box_{\mathfrak{F}}\blackdiamond_{\mathfrak{F}}(\vert \alpha\vert  \cup \vert \beta\vert )$ 
    \end{itemize}
    Proof:
    \begin{itemize}
        \item The ``$\supseteq$'' direction of the ``$\bot$'' clause follows from the fact that for any $\langle\Gamma, \Delta\rangle \in W$, $\bot \in \Delta$ (since $W \subseteq W^d_{k}$); while the ``$\supseteq$'' direction follows from (1) and the fact that $\Box_{\mathfrak{F}}\emptyset$ is the minimum in $FP_{\mathfrak{F}}$. 
        
        \item The ``$\supseteq$'' direction of the ``$\to$'' clause follows from Implication Lemma (i.e., Lemma \ref{five existence lemmas}.2) and $W^b_{k} \subseteq W$; while the ``$\subseteq$'' direction uses the fact that for any $\langle\Gamma, \Delta\rangle \in W$, $\!\to\!(\Gamma) \subseteq \Delta$ (since $W \subseteq W^d_{k}$).
        
        \item The ``$\land$'' clause follows from the rules ($\land$R) and ($\land$L) together with the fact that for any $\langle\Gamma, \Delta\rangle \in W$, $\Gamma$ is $\,\vdash$-closed (since $W \subseteq W^d_{k}$). 
        
        \item The ``$\supseteq$'' direction of the ``$\lor$'' clause follows from (1) and $\vert \alpha\vert \cup\vert \beta\vert  \subseteq \vert \alpha\lor\beta\vert $. 
        For the ``$\subseteq$'' direction: Suppose $\langle\Gamma_1, \Delta_1\rangle, \langle\Gamma_2, \Delta_2\rangle \in W$ with $\alpha \lor \beta \in \Gamma_1$ and $\Gamma_1 \cap \Delta_2=\emptyset$. So, $\alpha \lor \beta \notin \Delta_2$, then, since
        $\Delta_2$ is  $\,\vdash$-downward closed and closed under $\lor$ (as $W \subseteq W^d_{k}$), we have $Th_{\,\vdash_k}(\{\alpha\}) \cap \Delta_2 = \emptyset$ or $Th_{\,\vdash_k}(\{\beta\}) \cap \Delta_2 = \emptyset$. Then, by Extension Lemma (i.e., Lemma \ref{five existence lemmas}.1), we can find $\langle\Gamma_3, \Delta_3\rangle\in W^b_{k}\subseteq W$ such that $\Gamma_3\cap \Delta_2=\emptyset$ and $\Gamma_3$ contains either $\alpha$ or $\beta$.
    \end{itemize}

    Then, together with the definition of $\Vert\cdot\Vert$, we can easily prove by induction that for each $\alpha \in Form$, $\vert \alpha\vert  = \Vert\alpha\Vert$.
    
    \item
    \begin{enumerate}
        \item By Definition \ref{def: PsRe, PsSy}, we need to show that for each $\langle\Gamma, \Delta\rangle \in W$,  $\langle\Gamma, \Delta\rangle\in \Box_{\mathfrak{F}}\Diamond_{\mathfrak{F}}\Box_{\mathfrak{F}}\blackdiamond_{\mathfrak{F}}\{\langle\Gamma, \Delta\rangle\}$. Suppose $\langle\Gamma,\Delta\rangle R \langle\Gamma',\Delta'\rangle$. Since $W\subseteq W^d_{k}$,
        using Pseudo-symmetry Lemma (i.e., Lemma \ref{five existence lemmas}.5), there exists \(\langle \Phi, I_{k}^{\bot}(\Phi) \rangle\in W^b_{k}\subseteq W\) such that \(\Gamma \subseteq \Phi \) and $\langle\Gamma', \Delta'\rangle R \langle\Phi, I_{k}^{\bot}(\Phi)\rangle$. Since \(\Gamma \subseteq \Phi \), we conclude that $\langle\Phi, I_{k}^{\bot}(\Phi)\rangle \in \Box_{\mathfrak{F}}\blackdiamond_{\mathfrak{F}}\{\langle\Gamma, \Delta\rangle\}$.
        
        \item By Definition \ref{def: PsRe, PsSy}, we need to show that for each $\langle\Gamma, \Delta\rangle \in W$, $\langle\Gamma, \Delta\rangle\in \Box_{\mathfrak{F}}\emptyset\cup\Diamond_{\mathfrak{F}}\Box_{\mathfrak{F}}\blackdiamond_{\mathfrak{F}}\{\langle\Gamma, \Delta\rangle\}$. Suppose $\langle\Gamma, \Delta\rangle\notin \Box_{\mathfrak{F}}\emptyset$. Then, $\bot\notin\Gamma$ (for otherwise $\langle\Gamma, \Delta\rangle$ would have no successor). By Pseudo-reflexivity Lemma (i.e., Lemma \ref{five existence lemmas}.4), there exists \(\langle \Phi, I_{k}^{\bot}(\Phi) \rangle\in W^b_{k}\subseteq W\) such that \(\Gamma \subseteq \Phi \) and $\langle\Gamma, \Delta\rangle R \langle\Phi, I_{k}^{\bot}(\Phi)\rangle$. Since \(\Gamma \subseteq \Phi \), we conclude that $\langle\Phi, I_{k}^{\bot}(\Phi)\rangle \in \Box_{\mathfrak{F}}\blackdiamond_{\mathfrak{F}}\{\langle\Gamma, \Delta\rangle\}$. 
        
        \item Reflexivity follows easily from the fact that for each $\langle\Gamma, \Delta\rangle \in W$, $\Gamma\cap \Delta=\emptyset$ (since $W \subseteq W^d_{2}$).
    \end{enumerate}
    \end{enumerate}
\end{proof}

\begin{theorem}[Completeness]\label{completeness}
For each $k\in\{1,2\}$, $\Gamma\subseteq Form$, and $\varphi\in Form$,
    \[\Gamma\,\vDash_k\varphi \Rightarrow\Gamma\,\vdash_k\varphi.\]
\end{theorem}
\begin{proof}    
    Suppose $\Gamma\,\nvdash_k\varphi$. Then by Extension Lemma (i.e., Lemma \ref{five existence lemmas}.1),  there exists $\langle \Phi, I^{\bot}_k(\Phi) \rangle\in W^b_{k}$ such that $\Gamma \subseteq \Phi \text{ and }  \varphi\notin\Phi$.  Take any $W$ such that $W^b_{k} \subseteq W \subseteq W^d_{k}$, and let $\mathfrak{M}$ be the submodel of $\mathfrak{M}^d_{k}$ on $W$. Then, by Lemma \ref{Main lemma for the canonical model}, we know that $\mathfrak{M}\in\mathcal{D}_k$,  $\mathfrak{M},\langle \Phi, I^{\bot}_k(\Phi) \rangle\vDash\Gamma$, and $\mathfrak{M},\langle \Phi, I^{\bot}_k(\Phi) \rangle\nvDash\varphi$. Thus, $\Gamma\nvDash_k\varphi$.
\end{proof}

\appendix
\setcounter{definition}{0}
\setcounter{lemma}{0}

\renewcommand{\thedefinition}{A\arabic{definition}}
\renewcommand{\thelemma}{A\arabic{lemma}}

\section{}
In this appendix, we prove the completeness of $\,\vdash_{\textbf{K}}$ with respect to the class of all $\Box\blackdiamond$-models (cf. Theorem \ref{soundness and completeness of vdash_K}).

The approach is similar to that used in Section \ref{sec: completeness}. First, we define the $W^b$ and $W^d$ sets for $\,\vdash_{\textbf{K}}$:

\begin{definition}
    Define \(\mathfrak{M}^d_0 = \langle W^d_0, R^d_0, V^d_0 \rangle\) as follows:
    \begin{itemize}
        \item $W^d_0=\{\langle \Gamma, \Delta \rangle\in (\wp(Form))^2  \mid \Gamma$ is $\,\vdash_{\textbf{K}}$-closed,  $\{\bot\}\cup\!\to\!(\Gamma)\subseteq\Delta$, and $\Delta$ is $\,\vdash_{\textbf{K}}$-downward closed and closed under disjunction$\}.$

        \item \(R^d_0 = \{\langle \langle \Gamma_1, \Delta_1 \rangle, \langle \Gamma_2, \Delta_2 \rangle \rangle \in (W^d_0)^2 \mid \Gamma_1 \cap \Delta_2 = \emptyset\}\).
    
        \item \(V^d_0\) maps each $p\in PL$ to \( \{\langle \Gamma, \Delta \rangle \in W^d_0 \mid p \in \Gamma\}\).
    \end{itemize}
\end{definition}
Note that we no longer require that $\neg(\Gamma)\subseteq \Delta$ in the definition of $W^d_0$, since this requirement arose from the successor-seriality condition (cf. Lemma \ref{semantics: neg = to bot}), which is not presupposed on the semantics here.

\begin{definition}\label{definition of I_0(Gamma)}
Let
$I_0^{\bot}(\Gamma) = \{\varphi \mid$ there exists $\alpha_1,\dots,\alpha_n \in \Gamma$ and $\beta_1, \dots, \beta_n \notin  \Gamma$  such that $\varphi \vdash_{\textbf{K}} (\alpha_1 \rightarrow \beta_1) \lor \dots \lor (\alpha_n \rightarrow \beta_n) \lor \bot\}$.
Then define 
\begin{itemize}
    \item $W^b_{0}=\{\langle \Gamma,\, I_0^{\bot}(\Gamma) \rangle \mid  \Gamma\text{ is } \vdash_{\textbf{K}} \text{-closed and }\bot\notin\Gamma\}\, \cup\, \{\langle Form,\, \{\varphi \mid\varphi\vdash_{\textbf{K}} \gamma\}\rangle\}$.
\end{itemize}
\end{definition} 

Then we can check that: 
\begin{lemma}
    $W^b_{0}\subseteq W^d_{0}$.
\end{lemma}

Next, we prove the existence lemmas for $\,\vdash_{\textbf{K}}$ (similar to Lemma \ref{five existence lemmas}), except that we no longer need Pseudoreflexivity Lemma and Pseudosymmetry Lemma.
\begin{lemma}\label{existence lemmas for vdash_K}
    ~
\begin{enumerate}
    \item \textbf{Extension Lemma}\\
          If $ I\subseteq Form$ is closed under disjunction, and if $Th_{\,\vdash_{\textbf{K}}}(\Gamma)\cap I=\emptyset$, then there exists $\langle \Phi, \Psi \rangle\in W^b_{0}$ such that $\Gamma \subseteq \Phi \text{ and }  \Phi\cap I=\emptyset$.   
    \item \textbf{Implication Lemma}\\
          If $\Gamma\nvdash_{\textbf{K}}\alpha \to \beta$, then there is $\langle \Phi, I^{\bot}_0(\Phi) \rangle\in W^b_{0}$ such that $\Gamma \cap I^{\bot}_0(\Phi) = \emptyset \text{, } a \in \Phi \text{ and } \beta \notin \Phi$.
    \item \textbf{Fixpoint Lemma}\\
          If $\Gamma\nvdash_{\textbf{K}}\gamma $, then there exists  $\langle \Phi, \Psi \rangle\in W^b_{0}$ such that $\gamma\in\Psi$ and $\Gamma \cap \Psi = \emptyset$.
\end{enumerate}
\end{lemma}

\begin{proof}
    \begin{enumerate}
        \item If $I\neq\emptyset$, then take $\langle\Phi,\Psi\rangle=\langle Th_{\,\vdash_1}(\Gamma),\,I^{\bot}_0(Th_{\,\vdash_1}(\Gamma))\rangle$. Otherwise, take $\langle\Phi,\Psi\rangle=\langle Form,\, \{\varphi \mid\varphi\vdash_{\textbf{K}} \bot\}\rangle$.

        \item Assume $\Gamma\nvdash_{\textbf{K}}\alpha \to \beta$. Define the relation $\vsim$: $\alpha\vsim\beta\iff\Gamma\vdash_{\textbf{K}}\alpha\to\beta$. By Lemma \ref{three instances of vsim}, we know that $\,\vdash\,\subseteq\, \vsim$, and $\vsim$ satisfies (A), (Cut), ($\bot$), ($\land$L/R), ($\lor$R), and ($\Gamma\vdash_{\textbf{K}}\! f_{\bot}$).
        
        Let $\mathcal{S}=\{\Delta\subseteq Form \mid \Delta$ is $\vsim$-closed, $\alpha \in \Delta$, $\beta\notin\Delta\}$. 
        Since $\vsim$ satisfies (A), (Cut) and ($\land$L/R), and since $\Gamma\nvdash_{\textbf{K}}\alpha \to \beta$, we can prove that $Th_{\vsim}(\{\alpha\}) \in \mathcal{S}$. For any chain $H \subseteq \mathcal{S}$, clearly $\bigcup H \in \mathcal{S}$. 
        So, by Zorn's lemma, $\mathcal{S}$ has a maximal element, say $\Phi$, then by Lemma \ref{the encompassing lemma}, we have $\Gamma \cap I^{\bot}_0(\Phi) = \emptyset$. Moreover, since $\vsim$ satisfies ($\bot$), and since $\Phi$ is $\vsim$-closed and $\beta\notin\Phi$, we get $\bot\notin \Phi$. Thus, $\langle \Phi, I^{\bot}_0(\Phi) \rangle\in W^b_{0}$ is as required.

        \item Simply, $\langle Form,\, \{\varphi \mid\varphi\vdash_{\textbf{K}} \gamma\} \rangle\in W^b_{0}$ is as required.
    \end{enumerate}
\end{proof}

The proof of the following lemma is similar to that of Lemma \ref{Main lemma for the canonical model}. 
\begin{lemma}
Take any $W$ such that $W^b_{0} \subseteq W \subseteq W^d_{0}$, and let $\mathfrak{M} = \langle\mathfrak{F},V\rangle= \langle W,R,V\rangle$ be the submodel of $\mathfrak{M}^d_{0}$ on $W$. Define 
 $\vert \alpha\vert _W = \{\langle\Gamma, \Delta\rangle \in W \mid \alpha \in \Gamma\}$ for each formula $\alpha$. Then we have:
\begin{enumerate}[(1)]
    \item 
    For each $\alpha \in Form$, $\vert \alpha\vert _{W} \in FP_{\mathfrak{F}}$.
    \item \textbf{Truth Lemma:} For each $\alpha \in Form$, $\vert \alpha\vert _{W} = \Vert\alpha\Vert_{\mathfrak{M}}$.
\end{enumerate}
\end{lemma}

Then we obtain the completeness of $\,\vdash_{\textbf{K}}$:
\begin{theorem}[Completeness of $\,\vdash_{\textbf{K}}$]
For each $\Gamma\subseteq Form$ and $\varphi\in Form$,
    \[\Gamma\,\vDash_{\textbf{K}}\varphi \Rightarrow\Gamma\,\vdash_{\textbf{K}}\varphi.\]
\end{theorem}

%
%
%
\bibliographystyle{splncs04}
\bibliography{reference_2026.1}
\end{document}